%% file: main.tex
\documentclass[10pt]{article}

\usepackage{amsmath}
\usepackage{amssymb}
\usepackage{amsthm}
\usepackage[english]{babel}
\usepackage{dsfont}
\usepackage{float}
\usepackage{hyperref}
\usepackage{subcaption}
\usepackage{tikz}
\usepackage{url}

\usepackage{cleveref}
\usepackage[utf8]{inputenc}
\usepackage{pgfplots}

\usepackage[a4paper, total={30.5pc, 51pc}]{geometry}
\hypersetup{
	colorlinks=true,
	linkcolor=red,
	filecolor=magenta,      
	urlcolor=cyan,
	citecolor=green,
	linktoc=all,
}
\pgfplotsset{compat=1.18}

\crefformat{equation}{(#2#1#3)}

\theoremstyle{plain}
\newtheorem{theorem}{Theorem}[section]

\newtheorem{corollary}[theorem]{Corollary}
\newtheorem{lemma}[theorem]{Lemma}
\newtheorem{definition}[theorem]{Definition}

\theoremstyle{definition}
\newtheorem{remark}[theorem]{Remark}

\DeclareMathOperator*{\argmin}{arg\,min}
\DeclareMathOperator*{\Minimize}{\textsf{Minimize}}

\newcommand{\norm}[1]{\left\lVert#1\right\rVert}

\usepackage{mathtools}

\usepackage{tabularx}

\title{Resilient State Recovery using Prior Measurement Support Information}

\author{Yu Zheng\thanks{Department of Electrical and Computer Engineering, Florida State University, Tallahassee, FL 32310, USA
  (\url{yzheng6@fsu.edu}).}
\and Olugbenga Moses Anubi\thanks{Corresponding author. Department of Electrical and Computer Engineering, Florida State University, Tallahassee, FL 32310, USA.
  (\url{oanubi@fsu.edu}).}
\and Warren E. Dixon\thanks{Department of Mechanical and Aerospace Engineering, University of Florida, Gainesville FL, 32611-6250, USA
 (\url{wdixon@ufl.edu}). }
 }
\date{}

\begin{document}

\maketitle

\begin{abstract}
Resilient state recovery of cyber-physical systems has attracted much research attention due to the unique challenges posed by the tight coupling between communication, computation, and the underlying physics of such systems. By modeling attacks as additive adversary signals to a sparse subset of measurements, this resilient recovery problem can be formulated as an error correction problem. To achieve exact state recovery, most existing results require less than $50\%$ of the measurement nodes to be compromised, which limits the resiliency of the estimators. In this paper, we show that observer resiliency can be further improved by incorporating data-driven prior information. We provide an analytical bridge between the precision of prior information and the resiliency of the estimator. By quantifying the relationship between the estimation error of the weighted $\ell_1$ observer and the precision of the support prior. This quantified relationship provides guidance for the estimator's weight design to achieve optimal resiliency. Several numerical simulations and an application case study are presented to validate the theoretical claims.
\end{abstract}

\section{Notation and Preliminaries}\label{sec:notation}
\input{SICON/revised_sections/2-preliminary}

\section{Introduction}\label{sec:introduction}

\input{SICON/revised_sections/1-introduction}

\section{Model Development and Resilience Properties}\label{sec:Model}
\input{SICON/revised_sections/3-model}

\section{Vulnerability Analysis}\label{sec:vulnerability}
\input{SICON/revised_sections/4-vulnerability}

\section{Resilient State Recovery}\label{sec:Resilient_analysis}
\input{SICON/revised_sections/5-RESP}

\section{Numerical Simulation}\label{sec:simulation}

\input{SICON/revised_sections/6-simulation}

\section{Application Simulation}
\input{SICON/revised_sections/7-case_study}

\section{Conclusion}\label{sec:conclusion}
An enhanced resilient state estimation scheme is developed that uses prior information to improve the resilience of the weighted $\ell_1$ estimator. CSP and Row-RIP are properties introduced to measure the system's resilience under attacks. The traditional 2s-observability has also been shown to be one case of the proposed CSP. An optimization-based program is developed to generate a feasible attack against an $\ell_1$ estimator while the system's CSP is violated. The improvement of the weighted $\ell_1$ estimator's resilience is quantified with respect to the precision of the prior. The quantified relationship provides a theoretical guarantee of weight design for the weighted $\ell_1$ estimator. 

Through the analysis in this paper, we conclude that the precision of the prior is positively correlated to the resilience of the weighted $\ell_1$ estimator. The enhanced resilient estimation scheme is capable of maintaining resilient state recovery even when more than half of the sensors are attacked. The fundamental limitation of resilient recovery is obviated by incorporating data-driven prior into the physics-driven estimation design.

This paper has established a quantitative bridge that provides an in-road to many avenues for future research. For instance, given a target estimation performance and detection precision, how many sensors should be deployed in the system? The precision of a learning-based attack detector depends on the representative quality of the training dataset, so what is the quantitative relationship between the generalization of the generative attack dataset and the system's resiliency? Additionally, a recursive algorithm is expected for the weighted $\ell_1$ estimation scheme. For nonlinear systems, a similar quantitative bridge should be built for different linearization operators. While data-driven model predictive control (MPC) and moving-horizon estimation (MHE) have attracted much attention due to the difficulty of obtaining high-fidelity models, the resiliency of runtime data used for the data-driven models is not guaranteed. Therefore, a resilient data-driven MHE is expected, and the work presented in this paper could be utilized to achieve this goal.

\section*{Acknowledgments}
The authors acknowledge Professor Rodney Roberts at Florida State University for insightful discussions on the topic.

\bibliographystyle{unsrt}
\bibliography{references}

\end{document}

%% file: SICON/revised_sections/2-preliminary.tex
 ${\mathbb R}, {\mathbb R}^n, {\mathbb R}^{n \times m}$, $\mathbb{Z}$ denote the space of real numbers, real vectors of length $n$, real matrices of $n$ rows and $m$ columns, and the space of integers respectively. ${\mathbb R}_{+}$ denotes the space of positive real numbers. $\emptyset$ denotes the empty set. Normal-face lower-case letters $(e.g.\hspace{1mm} x \in {\mathbb R})$ are used to represent real scalars, bold-face lower-case letters $(e.g.\hspace{1mm} \mathbf{x} \in {\mathbb R}^n)$ represent vectors, while normal-face upper-case letters $(e.g.\hspace{1mm} X \in {\mathbb R}^{n \times m})$ represent matrices. $X^{\top}$ denotes the transpose of the matrix $X$. Let $\mathcal{T} \subseteq \{1,\cdots, n\}$, then $X_{\mathcal{T}} \in {\mathbb R}^{|\mathcal{T}| \times n}$ is the sub-matrix obtained by extracting the rows of $X \in {\mathbb R}^{m \times n}$ corresponding to the indices in $\mathcal{T}$. $\mathcal{T}^c$ denotes the complement of a set $\mathcal{T}$, and the universal set on which it is defined will be clear from the context. Also, $\mathcal{T}_1 \backslash \mathcal{T}_2 \triangleq \mathcal{T}_1 \cap \mathcal{T}_2^c$ denotes set difference between the sets $\mathcal{T}_1$ and $\mathcal{T}_2$. For a matrix $X\in \mathbb{R}^{m\times n}$, $\textbf{null}(X)$ and $\textbf{range}(X)$ denote its null space and range space respectively.
 
 For any $0<p\le\infty$, we use $\left\|\mathbf{x}\right\|_p \triangleq \left(\sum\limits_{i=1}^{n}{\lvert\mathbf{x}_i\rvert^p}\right)^{\frac{1}{p}}$ to denote the $p$-norm of the vector $\mathbf{x}\in\mathbb{R}^n$. A weighted $1$-norm of a vector $\mathbf{z} \in {\mathbb R}^{n}$ with the weight vector $\mathbf{w} \in \mathbb{R}^n$ is defined as: 
\begin{equation}
\|\mathbf{z}\|_{1, \mathbf{w}} \triangleq \sum_{i=1}^n \mathbf{w}_i|\mathbf{z}_i|.
\end{equation} 
We define a weight vector $\mathbf{w}(\mathcal{T},\omega) \in\mathbb{R}^m$, based on a set $\mathcal{T}\subset \{1,2,\cdots,m\}$, as:
\begin{equation} \label{equ:weight_definition}
    \mathbf{w}(\mathcal{T},\omega)_i =\left\{
     \begin{array}{lr}
     \omega &\quad\text{if}\hspace{0.2cm} i \in \mathcal{T} \\
     1 &\quad\text{otherwise}.
     \end{array}
     \right.
\end{equation}
where $0<\omega<1$ is the weight value.

The following mathematical preambles are used throughout the paper. The omitted proofs can be found in the standard references such as \cite{horn2012matrix}. For any vectors $\mathbf{x}, \mathbf{y}\in \mathbb{R}^n$ and real numbers $p,q\in \mathbb{R}_+$, the followings hold:
\begin{enumerate}
    \item The topological equivalence of vector norms:
        \begin{equation}
        \norm{\mathbf{x}}_q\le\norm{\mathbf{x}}_p\le n^{\left(\frac{1}{p}-\frac{1}{q}\right)}\norm{\mathbf{x}}_q,\hspace{1mm} \text{for} \hspace{1mm} 0<p\le q\le\infty.
        \end{equation}
        For example, given a vector $\mathbf{x} \in \mathbb{R}^n$,
        \begin{equation} \label{equ:topo_equi_norm}
        \|\mathbf{x}\|_\infty \le \|\mathbf{x}\|_2 \le \|\mathbf{x}\|_1 \le \sqrt{n}\|\mathbf{x}\|_2 \leq n\|\mathbf{x}\|_\infty.
        \end{equation}
    \item The triangle inequality:
        \begin{equation}\|\mathbf{x}+\mathbf{y}\|_p \leq \|\mathbf{x}\|_p + \|\mathbf{y}\|_p ,\hspace{1mm} \text{for} \hspace{1mm} 0<p\le \infty.
        \end{equation}
    \item The reverse triangle inequality:
        \begin{equation} \|\mathbf{x}-\mathbf{y}\|_p \geq \Big\lvert\|\mathbf{x}\|_p - \|\mathbf{y}\|_p\Big\rvert , \hspace{1mm} \text{for} \hspace{1mm} 0<p\le \infty.
        \end{equation} 
\end{enumerate}
\noindent The support of a vector $\mathbf{x} \in {\mathbb R}^{n}$ is a set of the indices of nonzero entries in $\mathbf{x}$, defined as
\begin{equation}
    \textsf{supp}(\mathbf{x}) \triangleq \{i\subseteq \{1,\dots, n\}\lvert\mathbf{x}_{i} \neq 0\},
\end{equation} 
and $\|\mathbf{x}\|_0 = |\textsf{supp}(\mathbf{x})|$. The vector $\mathbf{x}\in\mathbb{R}^n$ is said to be $k$-sparse if $\|\mathbf{x}\|_0 \le k$, and $\Sigma_{k}$ denotes the set of $k$-sparse vectors, given by
\begin{equation}\label{equ:sparsit_space}
    \Sigma_k \triangleq \left\{\mathbf{x}\in\mathbf{R}^n \bigg\lvert \hspace{2mm}|\textsf{supp}(\mathbf{x})| \le k\right\}.
\end{equation} 
Let $I_T \triangleq [i-T+1,i]$ denote a moving interval of fixed size $T$. When clear from context, the subscript $T$ will be dropped for clarity. Accordingly, $I-1 \triangleq [i-T,i-1]$ denotes the $T (>1)$ time interval from $i-T$ to $i-1$.
As a result, $\mathbf{x}_I = [\mathbf{x}_{i}^{\top}, \mathbf{x}_{i-1}^{\top}, \cdots, \mathbf{x}_{i-T+1}^{\top}]^{\top} \in \mathbb{R}^{Tn}$ is a flattened vector composed of vectors $\mathbf{x}_j$ in the descending order of $j\in I$. Given the supports $\mathcal{T}_j\subset \{1,2,\cdots,m\}, \forall j\in I $, the support sequence for the entire time interval $I$ is denoted by $\mathcal{T}_I = [\mathcal{T}_{i}^{\top}, \hspace{2mm} m+\mathcal{T}_{i-1}^{\top}, \cdots, (T-1)m+\mathcal{T}_{i-T+1}^{\top}]^{\top}$, where $m+\mathcal{T}_{i-1}^{\top}$ denotes adding $m$ to each element of $\mathcal{T}_{i-1}^{\top}$.

%% file: SICON/revised_sections/1-introduction.tex
The rapid development and wide application of cyber-physical systems (CPS) has benefited from the tight coupling and coordination among computation, communication, and physical components \cite{lee2016introduction}. Digital coupling via communication networks or physical coupling via grid systems provides wider observability, stronger controllability, and faster response speed among physical components. However, the closed-loop interaction between the cyber world and the physical world exposes the physical plant to cyberattacks \cite{pasqualetti2015control}. Small biases, injected through internet-of-thing sensors, could manipulate the physical processes maliciously \cite{zhang2020false, zhu2014attack}. Thus, for guaranteed resilient operation, it is necessary that true signals are recovered under cyber-activities and fault-induced anomalies. 

From a control system perspective, undesirable effects created by malicious intent could manifest as noise, disturbance, uncertainty, or any combination. However, attack signals are more challenging since they are possibly unbounded, or bounded with unknown bounds, and time-varying according to the attacker's policy \cite{sui2020vulnerability}. These special features of attack signals result in the typical robust control and estimation laws \cite{wang2022robust}, where bound information is required, being ineffective \cite{weerakkody2019resilient}. Assuming attackers don't have infinite resources, the attack vectors $\mathbf{e} \in \mathbb{R}^m$ are sparse with bounded support $\|\mathbf{e}\|_0 \leq s < m$ \cite{fawzi2014secure}. Based on this sparsity assumption on the set of attacked nodes, the majority of state recovery results are formulated as a sparse recovery problem closely related to the classical error correction problem \cite{candes2005decoding, fawzi2014secure}. 

Consider a measurement model of the form
\begin{equation}
    \mathbf{y} = H\mathbf{x}+\mathbf{e},
\end{equation}
where $H \in \mathbb{R}^{m\times n}$ ($n<m$) is the measurement matrix, and $\mathbf{y},\mathbf{e}, \mathbf{x}$ are the sensor measurements, attack injection, and internal state vectors respectively. Given a coding matrix $F\in \mathbb{R}^{n\times m}$ satisfying $FH=0$, the classical error correction problem is to recover a sparse vector $\mathbf{e}$ for which $\mathbf{y}_c \triangleq F\mathbf{y} = F\mathbf{e}$. This is usually solved by 
\begin{equation}\label{equ:0-norm}
    \underset{\mathbf{e}}{\textsf{Minimize}}: \|\mathbf{e}\|_0 \hspace{2mm} \textsf{Subject to}: \mathbf{y}_c = F\mathbf{e}.
\end{equation}
Any sparse error $\|\mathbf{e}\| < s \hspace{2mm}(s<\frac{m}{2})$ could be recovered uniquely through the above optimization program if \cite{vidyasagar2019introduction}
\begin{equation}\label{equ:uniqueness_condition}
    \textsf{null}(F) \cap \Sigma_{2s} = \emptyset ,
\end{equation}
where $\Sigma_{2s}\subseteq\mathbb{R}^m$ denotes the set of all $2s$-sparse vectors which is a subset of the $m$-dimensional vector space. The minimization problem in \eqref{equ:0-norm} is an NP-hard problem but can be relaxed to its convex neighbor:
\begin{equation}\label{equ:1-norm}
    \underset{\mathbf{e}}{\textsf{Minimize}}: \|\mathbf{e}\|_1 \hspace{2mm} \textsf{Subject to}: \mathbf{y}_c = F\mathbf{e}.
\end{equation}
The solution of \eqref{equ:1-norm} is the unique solution of \eqref{equ:0-norm} if the coding matrix $F$ satisfies the null space property (NSP) \cite{candes2006robust} or s-restricted isometric property (s-RIP) \cite{candes2005decoding}. The authors in \cite{pajic2016attack} also studied the robustness of $\ell_1$ estimators while attacks and noise both exist. Alternative approaches employ $\ell_2$ decoders as the convex relaxation of the program in \eqref{equ:0-norm}. These methods utilize diverse sparse projectors, such as an event-trigger or switch mechanisms \cite{an2017secure, shoukry2015event} and constrained sensor fusion \cite{chen2018resilient, nakahira2018attack, shoukry2017secure}. These studies have also been extended to networked control systems \cite{ mitra2019byzantine, sargolzaei2021lyapunov, sargolzaei2019detection, sargolzaei2022secure}. However, both relaxed convex programs admit the uniqueness condition in \eqref{equ:uniqueness_condition}, regardless of whether they are implemented in centralized or distributed systems. This condition can be cast as a $2s$-observerability assumption for linear systems, requiring more than half of the sensors to be safe \cite{fawzi2014secure,shoukry2015event}. In other words, $100\%$ sensor redundancy is required. This assumption poses a major limitation on CPS's resiliency in practice. 

The fundamental limitation in \eqref{equ:uniqueness_condition} can be relaxed by incorporating prior information. The authors in \cite{friedlander2011recovering} studied a weighted $\ell_1$ decoder with support prior for compressed sensing problems. The authors in \cite{vaswani2010modified} studied a modified compressed sensing program to recursively reduce the state estimation error along with convergence of the support prior estimation. In \cite{khazraei2022attack}, the authors found the resilience of the state estimator can be improved by incorporating intermittent data authentication. The authors in \cite{shinohara2019resilient} also found that utilizing the state prior could also increase the resilience of the state estimator. However, analysis is needed for why such methodology can improve resiliency and how much resiliency can be improved.

\textbf{Contributions}: In this paper, we studied an enhanced resilient estimation problem with reinforcement of data-driven prior information. Existing resilient estimation methods predominantly rely on the 2s-observability assumption \cite{an2017secure, chen2018resilient, fawzi2014secure, he2021distributed, lee2020fully, mitra2019byzantine, nakahira2018attack, shoukry2015event, shoukry2017secure}. Our work seeks to relax this assumption, thereby further improving the resiliency of CPS. To achieve this goal, we derive the column space property (CSP) and row-restricted isometry property (Row-RIP), which enables us to adjust the strictness of the observability conditions in relation to the estimation error bound. Our demonstrations indicate that resiliency can be significantly enhanced when the accuracy of the prior information exceeds that of a random coin flip. Furthermore, while previous studies have explored resilient estimation with additional prior information \cite{khazraei2021learning, shinohara2019resilient}, they do not fully explain why and how this information improves resiliency, nor do they identify the necessary conditions for its effectiveness. Our paper addresses this gap by establishing a connection between the precision of data-driven prior information and the resulting estimation error bound. These quantitative connections provide a systematic guide for the design of resilient estimators leveraging prior information,  as illustrated by the design of the weighted $\ell_1$ estimator presented in this paper.
 
 The rest of the paper is organized as follows. Notations and preliminary results are provided in Section~\ref{sec:notation}. In Section~\ref{sec:Model}, we introduce CPS and Row-RIP, motivated by the uniqueness condition of state recovery. In Section~\ref{sec:vulnerability}, we conduct a rigorous vulnerability analysis and show that violating the CSP implies the existence of a successful attack design against a resilient $\ell_1$ estimator. We derive the estimation errors bound for the $\ell_1$ decoder without prior and weighted $\ell_1$ decoder with prior in Section~\ref{sec:Resilient_analysis} and show that inclusion of prior information significantly lowers the resulting error bound. In Section~\ref{sec:simulation}, a numerical simulation quantifies the influence of the precision of the prior information on system resilience. Concluding remarks are in Section~\ref{sec:conclusion}.

%% file: SICON/revised_sections/3-model.tex
In this section,  we present the model and relevant properties of a CPS. Moreover, instead of assuming a $2s$-observability as done in results such as \cite{fawzi2014secure, shoukry2015event}, we introduce CSP and Row-RIP to guarantee the uniqueness of state recovery under sparse attacks.

\subsection{System Model}\label{System Model with Bad Data Detection}
Since the physical processes interact with the cyber components that operate in a discrete-time sequence, we consider the following discrete model of a CPS under attack and bounded noise: 
\begin{equation} \label{equ:system_model}
 \begin{aligned}
     \mathbf{x}_{i+1}&=A\mathbf{x}_i\\
    \mathbf{y}_i&=C\mathbf{x}_i+\mathbf{e}_i,
 \end{aligned}
\end{equation}
where $\mathbf{x}_i \in {\mathbb R}^{n}, \mathbf{y}_i \in {\mathbb R}^{m}$ are the state and measurement vector at the time instance $i$ respectively, $\mathbf{e}_i \in \mathbb{R}^m$ is a vector containing both noise and attack with attack support $\mathcal{T} \subset \{1,2,\cdots,m\}$ ($|\mathcal{T}|<k$). The noise portion of $\mathbf{e}$ is assumed to be bounded
\begin{equation}\label{equ:noise_bound}
\sum\limits_{i\in\mathcal{T}^c}|\mathbf{e}_i|<\varepsilon,
\end{equation} 
where $\varepsilon \in \mathbb{R}_{+}$ is the noise level. The known control inputs are neglected in \eqref{equ:system_model} since they are irrelevant to state estimation problem \cite{fawzi2014secure,zheng2021Resilient} or the reader could view \eqref{equ:system_model} as a stable closed-loop dynamical model. By iterating \eqref{equ:system_model} $T$ steps backwards, the measurement model in $T$-length observation interval is given by
\begin{equation}\label{equ:measurement_model}
    \mathbf{y}_I = H\mathbf{x}_{i-T+1}+\mathbf{e}_I,
\end{equation}
where $\mathbf{y}_I, \mathbf{e}_I \in \mathbb{R}^{Tm}$ are the sequence of measurements and attacks in the interval $I$ respectively, $\mathbf{x}_{i-T+1}\in \mathbb{R}^n$ is the state vector at time $i-T+1$, and $H = \begin{bmatrix}  (C A^{T-1})^{\top} & \cdots & (CA)^{\top} & C^{\top} \end{bmatrix}^{\top}$ is assumed to have full column rank. 

\noindent Due to the sparsity of the attack portion in $\mathbf{e}_I$, an $\ell_1$ decoder is typically used to reconstruct the state from the corrupted measurements \cite{anubi2019enhanced, candes2005decoding, fawzi2014secure}.
\begin{definition}[\textbf{Decoder of Horizon $T$}]\label{Def:decoder}
Given an observation matrix $H$, the decoder $\mathcal{D}$ maps the observation sequence $\mathbf{y}_I$ to a state estimate $\hat{\mathbf{x}}_{i-T+1}$. The state estimate is
\begin{equation}\label{equ:decoder}
 \begin{aligned}
    \hat{\mathbf{x}}_{i-T+1} &= \mathcal{D}(\mathbf{y}_I|H) \\
    & \triangleq \argmin_{\mathbf{x} \in \mathbb{R}^n} \|\mathbf{y}_I - H \mathbf{x}\|_1,
\end{aligned}
\end{equation}
where $\mathbf{y}_I$ is the observation history on $I$ and $\hat{\mathbf{x}}_{i-T+1}\in\mathbb{R}^n$ is the resulting estimated initial state vector.
\end{definition}
The $\ell_1$ decoder defined in \eqref{equ:decoder} can be efficiently solved by linear programming \cite{candes2006robust, candes2005decoding}, in a recursive manner \cite{needell2009noisy}, or through distributed solvers \cite{lee2020fully}. 
\begin{definition}[\textbf{Resilient Recovery}]
The decoder $\hat{\mathbf{x}}_{i-T+1} = \mathcal{D}(\mathbf{y}_I|H)$ in Definition~\ref{Def:decoder} is said to resiliently recover the true state $\mathbf{x}_{i-T+1}$ if, given an error tolerance $\epsilon$,
\begin{align}\label{equ:resilient_recovery}
    \|\hat{\mathbf{x}}_{i-T+1} - \mathbf{x}_{i-T+1} \| \leq \epsilon.
\end{align}
If \eqref{equ:resilient_recovery} holds, then the decoder $\mathcal{D}$ is a resilient decoder with error tolerance $\epsilon$.
\end{definition}

\noindent To achieve resilient recovery through the estimator in \eqref{equ:decoder}, the literature commonly assumes a condition, known as 2s-observability \cite{fawzi2014secure,lee2015secure}, to ensure the unique existence of a state despite the presence of attack signals. However, this assumption imposes a significant constraint on system resiliency, because it requires $100\%$ sensor redundancy relative to the attacker's capabilities. In the remainder of this work, we examine the original conditions that are sufficient for resilient recovery and explore the potential to relax the 2s-observability assumption, thereby enhancing the resiliency of CPS.

\subsection{CSP and Row-RIP}
We begin with the result that guarantees unique state recovery from attacked measurements. Without loss of generality, the attack portion and noise portion in the error term $\mathbf{e}$ can be seen as orthogonal to each other since the noise on the attack portion indexed by $\mathcal{T}$ can be absorbed into the attack itself. In addition, uniqueness in state recovery is a hallmark of sparse recovery problems when countering the effects of these attacks, as demonstrated in \cite{candes2008restricted, candes2005decoding}. Consequently, we focus our discussion on unique state recovery using the noise-free version of the measurement model in \eqref{equ:measurement_model}.

\begin{theorem}[\textbf{Uniqueness of State Recovery}] \label{Thm:uniquness}
Given an observation matrix $H \in \mathbb{R}^{Tm\times n}$ subject to $k$-sparse attacks on $T$-length time interval $I$ (i.e. $|\mathcal{T}_j| \leq k$ for all $j\in I$). If, for any nonzero $\mathbf{h} \in \textsf{range}(H)$, it is true that
\begin{equation}\label{equ:uniqueness}
\begin{aligned}
        \|\mathbf{h}_{\mathcal{S}}\|_1 < \|\mathbf{h}_{\mathcal{S}^c}\|_1, \hspace{1mm} \text{for all} \hspace{1mm} \mathcal{S} \subset \{1,2,\cdots,Tm\}\hspace{1mm}
    \text{with} \hspace{1mm} |\mathcal{S}| \leq Tk,
\end{aligned}
\end{equation}
then for each attacked measurement $\mathbf{y}_I \in \mathbb{R}^{Tm}$, there exists a unique state vector $\hat{\mathbf{x}} \in \mathbb{R}^{n}$ and an attack vector $\hat{\mathbf{e}}_I$ satisfying \eqref{equ:measurement_model}.
\end{theorem}
\begin{proof}
Let $(\mathbf{z}_1,\mathbf{e}_1), (\mathbf{z}_2,\mathbf{e}_2) \in \mathbb{R}^n \times \Sigma_{Tk}$ be two different candidates satisfying \eqref{equ:measurement_model} for a given $\mathbf{y}_I$. That is: $ \mathbf{y}_I = H\mathbf{z}_1+\mathbf{e}_1=H\mathbf{z}_2+\mathbf{e}_2$. Then $ H(\mathbf{z}_1-\mathbf{z}_2) = \mathbf{e}_2-\mathbf{e}_1$. This is equivalent to $H(\mathbf{z}_1-\mathbf{z}_2) \in \Sigma_{2Tk}$ for some $\mathbf{z}_1,\mathbf{z}_2 \in \mathbb{R}^n$. Thus, since $H$ is full-ranked, the uniqueness condition $\mathbf{z}_1 = \mathbf{z}_2$ holds if and only if $ \textsf{range}(H) \cap \Sigma_{2Tk} = \{\mathbf{0}\}.
$ Given a nonzero $\mathbf{h} \in \textsf{range}(H)$, this condition is equivalent to $\|\mathbf{h}\|_0>2Tk$. 

To prove $\|\mathbf{h}\|_0>2Tk$ can be implied by the condition in \eqref{equ:uniqueness}, we suppose, for the sake of contradiction, that $\|\mathbf{h}\|_0\leq 2Tk$. Choose $\overline{\mathcal{S}}\subseteq \{1,2,\cdots,Tm\}$, with $|\overline{\mathcal{S}}|=Tk$, to be the indices of the largest components of $\mathbf{h}$ in absolute value. Then, it must be that $\|\mathbf{h}_{\overline{\mathcal{S}}}\|_1 \geq \|\mathbf{h}_{\overline{\mathcal{S}}^c}\|_1,
$ which is a contradiction. Thus, \eqref{equ:uniqueness} implies that $\|\mathbf{h}\|_0>2Tk$.
\end{proof}
\begin{remark}\label{rmk:euqal_2s_observability}
According to the proof of Theorem~\ref{Thm:uniquness}, the inequality in \eqref{equ:uniqueness} implies $\|\mathbf{h}\|_0>2Tk$. This result is equivalent to $m>2k$ which is well-known and used in related work (cf. \cite{fawzi2014secure, shoukry2015event}).
\end{remark}

\begin{remark}
Using the compactness of the intersection of $\textsf{range}(H)$ and the unit norm ball, the inequality in \eqref{equ:uniqueness} is equivalent to the existence of $\beta \in (0,1)$ such that $ \|\mathbf{h}_{\mathcal{S}}\|_1 \leq \beta\|\mathbf{h}_{\mathcal{S}^c}\|_1$ holds for all $\mathbf{h} \in \textsf{range}(H)$, $\mathcal{S} \subset \{1,2,\cdots,Tm\}, |\mathcal{S}| \leq Tk$. 
\end{remark}

\begin{definition}[\textbf{Column space property (CSP)}]\label{def:CSP}
A matrix $H \in \mathbb{R}^{m\times n}$ is said to have a CSP of order $s<m$ (denoted as $H\triangleright \textsf{CSP}(s)$) with a parameter $\beta$ if, for every $\mathbf{h}\in \textsf{range}(H)$, 
\begin{equation}\label{equ:CSP}
\begin{aligned}
\|\mathbf{h}_{\mathcal{S}}\|_1 \leq \beta\|\mathbf{h}_{\mathcal{S}^c}\|_1, \hspace{1mm} \text{for all} \hspace{1mm} \mathcal{S} \subset \{1,2,\cdots,m\} \hspace{1mm} \text{with} \hspace{1mm} |\mathcal{S}|\leq s.
\end{aligned}
\end{equation}
\end{definition}
Definition~\ref{def:CSP} is similar to the well-known Null Space Property (NSP) \cite{candes2006robust} but is defined on the range space instead. The name CSP is used to avoid confusion with a recently defined Range Space Property (RSP) \cite{zhao2014equivalence}. Next, we consider a restricted isometry property on the row space of a matrix and establish a relationship with CSP.
 
 \begin{definition}[\textbf{Row-Restricted Isometry Property (Row-RIP)}]\label{def:ros-rip}
 A matrix $H \in \mathbb{R}^{m\times n}$ is said to have a Row-RIP of order $k$ (denoted as $H \triangleright \textsf{rRIP}(k)$) if there exists a constant $\delta \in (0,1)$ such that
 \begin{equation}\label{eqn:rRIP}
 (1-\delta)\|\mathbf{x}\|_2^2 \leq \|H_{\mathcal{T}}\mathbf{x}\|_2^2 \leq (1+\delta)\|\mathbf{x}\|_2^2,
 \end{equation}
 for all $\mathbf{x} \in \mathbb{R}^n$ and $\mathcal{T} \subset \{1,2,\cdots, m\}$ with $|\mathcal{T}| \leq k$. 
 Furthermore, the row-RIP constant of order $k$ is defined as
  \begin{equation}
 \delta_k \triangleq \inf_{\delta}\left\{\delta\mid \eqref{eqn:rRIP} \text{ holds }\right\}.
 \end{equation}
 \end{definition}
 \begin{lemma}[\textbf{Relating CSP and Row-RIP}]\label{Lem:rRIP_CSP}
Given $H\in \mathbb{R}^{m\times n}$, with $H \triangleright \textsf{rRIP}(ak)$ for some $a>1$ and $k\leq \frac{m}{2a+1}$. If
\begin{equation}\label{eqn:rRIP_CSP} 
\begin{aligned}
   \delta_k + a\delta_{ak} < a-1,
\end{aligned}
\end{equation}
then $H \triangleright\textsf{CSP}(k)$.
\end{lemma}
\begin{proof}
Choose an arbitrary set $\mathcal{T} \in \{1,2,\cdots,m\}$ with $|\mathcal{T}| \leq k \leq m/(2a+1)$, and an arbitrary vector $\mathbf{h} \in \textsf{range}(H)$. Then using the topology equivalence of norms in \eqref{equ:topo_equi_norm} and the inequality in \eqref{eqn:rRIP}, it follows that
\begin{equation} \label{equ:Lemma_equ1}
    \begin{aligned}
\|\mathbf{h}_{\mathcal{T}}\|_1 \leq \sqrt{|\mathcal{T}|}\|\mathbf{h}_{\mathcal{T}}\|_2\leq \sqrt{k(1+\delta_k})\|\mathbf{x}\|_2,
\end{aligned} 
\end{equation}
with $\mathbf{h} = H\mathbf{x}$. Also, let $\mathcal{S} \subset \mathcal{T}^c$, with $|\mathcal{S}|=ak$, be a subset of $\mathcal{T}^c$ containing the smallest entries, in absolute values, of $\mathbf{h}_{\mathcal{T}^c}$. Such set $\mathcal{S}$ exists because $|\mathcal{S}|+|\mathcal{T}|\leq ak+k\leq m$. Then, the largest element in $\mathbf{h}_{\mathcal{S}}$ satisfies
\begin{align*}
    \left|\mathcal{T}^c\backslash \mathcal{S}\right|\max_{i\in\mathcal{S}}\left|\mathbf{h}_i\right|\le\sum_{i\in\mathcal{T}^c\backslash\mathcal{S}}\left|\mathbf{h}_i\right|\le\sum_{i\in\mathcal{T}^c}\left|\mathbf{h}_i\right|,
\end{align*}
which implies that
\begin{equation} \label{equ:max_min_set}
  (m-|\mathcal{T}|-|\mathcal{S}|)\max_{i\in\mathcal{S}} |\mathbf{h}_i| \leq \|\mathbf{h}_{\mathcal{T}^c}\|_1.
\end{equation}
The inequality in \eqref{equ:max_min_set} follows from the topological equivalence of norms in \eqref{equ:topo_equi_norm}. Besides, since $k\le\frac{m}{2a+1}$ and $|\mathcal{S}| = ak$, it follows that $(2a+1)k\leq m \Rightarrow (2a+1)k-k-ak \leq m -|\mathcal{T}| -|\mathcal{S}| \Rightarrow m -|\mathcal{T}| -|\mathcal{S}|\geq ak $. Thus, it follows from \eqref{equ:max_min_set} that 
\begin{equation}\label{equ:Lemma23_proof_1}
      \max_{i\in\mathcal{S}} |\mathbf{h}_i| \leq \frac{\|\mathbf{h}_{\mathcal{T}^c}\|_1}{ak} \Rightarrow \|\mathbf{h}_{\mathcal{S}}\|_\infty \leq \frac{\|\mathbf{h}_{\mathcal{T}^c}\|_1}{ak}.
\end{equation}
By considering the topological equivalence of the $\infty$-norm and $2$-norm in \eqref{equ:topo_equi_norm}, \eqref{equ:Lemma23_proof_1} implies
\begin{equation}
    \sqrt{ak}\|\mathbf{h}_{\mathcal{S}}\|_2 \leq \|\mathbf{h}_{\mathcal{T}^c}\|_1.
\end{equation}
Since $H \triangleright \textsf{rRIP}(ak)$, it follows that
\begin{equation*}
    \|\mathbf{h}_{\mathcal{T}^c}\|_1 \geq {\sqrt{ak}} \|\mathbf{h}_{\mathcal{S}}\|_2 \geq \sqrt{ak(1-\delta_{ak})}\|\mathbf{x}\|_2,
\end{equation*}
which implies 
\begin{equation}\label{equ:Lemma_equ2}
    \|\mathbf{x}\|_2 \leq \frac{1}{\sqrt{ak(1-\delta_{ak})}} \|\mathbf{h}_{\mathcal{T}^c}\|_1.
\end{equation}
Combining \eqref{equ:Lemma_equ1} and \eqref{equ:Lemma_equ2} yields 
\begin{equation*} \|\mathbf{h}_{\mathcal{T}}\|_1 \leq \sqrt{\frac{1+\delta_{k}}{a(1-\delta_{ak})}}  \|\mathbf{h}_{\mathcal{T}^c}\|_1.
\end{equation*}
From \eqref{eqn:rRIP_CSP}, it follows that $0<\sqrt{\frac{1+\delta_{k}}{a(1-\delta_{ak})}}<1$. Thus, $H \triangleright \textsf{CSP}\left(k\right)$.
\end{proof}

%% file: SICON/revised_sections/4-vulnerability.tex
In this section, we discuss under what conditions the system is vulnerable to attacks in a conservative scenario. We make the following assumptions about the attackers' capability which are widely admitted in literature (cf. \cite{hu2018state, pasqualetti2015control} and references within):
 \begin{enumerate}
     \item The attacker has full model knowledge of \eqref{equ:system_model};
     \item  The attacker can inject arbitrary bias at the compromised sensors $\mathcal{T}_i \subset \{1,\cdots,m\}$;
     \item The resource of the attacker is limited, so the number of simultaneously attacked sensors is bounded such that $|\mathcal{T}_i| \leq k$ (equivalently, $\mathbf{e}_i \in \Sigma_k$). 
 \end{enumerate}

\noindent These assumptions highlight the insufficiency of conventional bounded-disturbance-based robust estimation approaches since the attacks may not be bounded. Furthermore, it is more pragmatic to assume that attackers aim to maximize the perturbation of estimated states without triggering an alarm from the monitoring system.

\begin{definition}[\textbf{Successful attack} \cite{ chen2022geometrical, hu2018state}]\label{Def:Successful_FDIA}
Consider the CPS in \eqref{equ:system_model} and the corresponding measurement model in \eqref{equ:measurement_model} with a state estimator $\mathcal{D}: \mathbb{R}^{Tm} \rightarrow \mathbb{R}^n$. Given an attack support $\mathcal{T}$ with $|\mathcal{T}_j|\leq k$ for all $j\in I$ and a positive integer $k<m$, the attack sequence $\mathbf{e}_I\in\mathbb{R}^{Tm}$ is said to be $(\epsilon, \alpha)$-successful if
\begin{equation}\label{equ:def_success_FDIA}
     \|\mathbf{x}^{\star}-\mathcal{D}(\mathbf{y}_I)\|_2 \geq \alpha \hspace{2mm} \text{and} \hspace{2mm}  \|\mathbf{y}_I-H\mathcal{D}(\mathbf{y}_I)\|_2 \leq \epsilon,
\end{equation}
where $\mathbf{y}_I = \mathbf{y}^{\star}_I+\mathbf{e}_I$ with $\mathbf{y}^{\star}_I\in\mathbb{R}^{Tm}$ being the true measurement vector, and $\mathbf{x}^{\star}=\mathcal{D}(\mathbf{y}^{\star}_I)$ is the observed state vector without attacks.
\end{definition}
The first inequality in \eqref{equ:def_success_FDIA} quantifies the attack effectiveness using a threshold value of $\alpha$ while the second inequality quantifies the stealthiness using a threshold value of $\epsilon$. Given a support sequence $\mathcal{T}$ with $|\mathcal{T}_i| \leq k$ for all $j\in I$. Let $\mathbf{x}_e$ be an optimal solution of
\begin{equation}\label{equ:FDIA_optimization}
\begin{aligned}
    \textsf{Maxmize}&: \hspace{0.2cm} \|H_{\mathcal{T}} \mathbf{x}\|_{1}, \\
    \textsf{Subject to}&: \hspace{0.2cm} \| H_{\mathcal{T}^c}\mathbf{x}\|_{1} \leq \epsilon.
\end{aligned}
\end{equation}
Consider the attack defined accordingly as
\begin{equation}\label{equ:FDIA}
    \mathbf{e}_{I_\mathcal{T}} = H_{\mathcal{T}} \mathbf{x}_e, \hspace{2mm} \mathbf{e}_{I_{\mathcal{T}^c}} = \mathbf{0}.
\end{equation}
The next result gives the conditions under which the attack defined above is $(\epsilon,\alpha)$-successful for the given attack support $\mathcal{T}$.
\begin{remark}
    It is implicitly assumed in \eqref{equ:FDIA} that $H_{\mathcal{T}^c}$ is full ranked. Otherwise, the constraint $ \| H_{\mathcal{T}^c}\mathbf{x}\|_{1} \leq \epsilon$ would be satisfied by any $\mathbf{x}_e$ in the null space of $H_{\mathcal{T}^c}$, thereby making the program in \eqref{equ:FDIA_optimization} ill-defined.
\end{remark}
\begin{theorem}[\textbf{$(\epsilon,\alpha)$-successful attack}]\label{Thm:FDIA}
Given an attack support $\mathcal{T}$, with $|\mathcal{T}|>\frac{\underline{\sigma}_{\mathcal{T}^c}^2}{\overline{\sigma}_{\mathcal{T}}^2}$, such that $H_{\mathcal{T}^c}$ has full rank. If the subset
\begin{equation}\label{equ:violation_NSP}
 \mathcal{S}_{H}(\mathcal{T}) \triangleq \{ \mathbf{w} \in \textsf{range}(H) \big| \|\mathbf{w}_{\mathcal{T}}\|_1 > \|\mathbf{w}_{\mathcal{T}^c}\|_1 \},
\end{equation}
is non-empty, then the attack in \eqref{equ:FDIA} is $(\epsilon,\alpha)$-successful for all 
\begin{equation}\label{equ:alpha_L1_FDIA}
    \alpha \leq \frac{\sigma_1-1}{\sqrt{|\mathcal{T}|}\overline{\sigma}_{\mathcal{T}} - \underline{\sigma}_{\mathcal{T}^c}}\epsilon,
\end{equation}
where $\overline{\sigma}_{\mathcal{T}}$ and $\underline{\sigma}_{\mathcal{T}^c}$ are the largest and the smallest nonzero singular values of $H_{\mathcal{T}}$ and $H_{\mathcal{T}^c}$ respectively, and
$ \sigma_1 = \displaystyle{\max_{\mathbf{v}\in \mathbb{R}^n \backslash \{\mathbf{0}\} }} \frac{ \|H_{\mathcal{T}}\mathbf{v}\|_1}{\|H_{\mathcal{T}^c}\mathbf{v}\|_1} .
$
\end{theorem}

\begin{proof}
First, we show that $ \|\mathbf{x}^{\star}-\mathcal{D}(\mathbf{y}_I)\|_2 \geq \alpha$. For the optimization problem in \eqref{equ:FDIA_optimization}, let $\mathbf{x} = \alpha \mathbf{v}$, where $\mathbf{v} \in \mathbb{R}^n$ is an arbitrary vector, then
\begin{equation}\label{equ:FDIA_Proof_1}
    \begin{aligned}
    \|H_{\mathcal{T}}\mathbf{x}_{e}\|_1 &\geq \max_{|\alpha|\|H_{\mathcal{T}^c}\mathbf{v}\|_1\leq \epsilon} \left(|\alpha|\|H_{\mathcal{T}}\mathbf{v}\|_1\right) \\
&\geq \max_{\mathbf{v}\in \mathbb{R}^n \backslash \{\mathbf{0}\} } \frac{ \|H_{\mathcal{T}}\mathbf{v}\|_1}{\|H_{\mathcal{T}^c}\mathbf{v}\|_1}\epsilon = \sigma_1 \epsilon
\end{aligned} 
\end{equation}
Let $P$ be an appropriate permutation matrix satisfying $H = P\begin{bmatrix}H_{\mathcal{T}} \\ H_{\mathcal{T}^c}\end{bmatrix}$. Then, according to \eqref{equ:FDIA}, $\mathbf{e}_I = P\begin{bmatrix}H_{\mathcal{T}} \\ 0\end{bmatrix}\mathbf{x}_e$. Consider a projection,
\begin{equation}\label{equ:definition_e_perp}
    \mathbf{x}_e^{\perp} \triangleq \argmin_{\mathbf{z}}\left\|H\mathbf{z}-P\begin{bmatrix}H_{\mathcal{T}} \\ 0\end{bmatrix}\mathbf{x}_e\right\|_1,
\end{equation}
it follows that 
\begin{equation*}
\begin{aligned}
    \left\|H\mathbf{x}_e^{\perp} - P\begin{bmatrix}H_{\mathcal{T}} \\ 0\end{bmatrix}\mathbf{x}_e\right\|_1 &\leq \left\|H\mathbf{x}_e - P\begin{bmatrix}H_{\mathcal{T}} \\ 0\end{bmatrix}\mathbf{x}_e\right\|_1 \\
&\leq \left\|H_{\mathcal{T}^c}\mathbf{x}_e \right\|_1 \leq \epsilon.
\end{aligned}
\end{equation*}
Then, by decomposing $1$-norm on the disjoint sets $\mathcal{T}$ and $\mathcal{T}^c$, we have that
\begin{equation*}
\begin{aligned}
     \left\|H\mathbf{x}_e^{\perp} - P\begin{bmatrix}H_{\mathcal{T}} \\ 0\end{bmatrix}\mathbf{x}_e\right\|_1 &= \left\|P\begin{bmatrix}H_{\mathcal{T}}(\mathbf{x}_e^{\perp} - \mathbf{x}_e) \\ H_{\mathcal{T}^c}\mathbf{x}_e^{\perp} \end{bmatrix} \right\|_1 \\
     &= \left\|H_{\mathcal{T}}(\mathbf{x}_e^{\perp}-\mathbf{x}_e)\right\|_1 +\left\|H_{\mathcal{T}^c}\mathbf{x}_e^{\perp}\right\|_1.
\end{aligned}
\end{equation*}
Thus, it follows that
\begin{equation}\label{equ:optimal_projection}
 \left\|H_{\mathcal{T}}(\mathbf{x}_e^{\perp}-\mathbf{x}_e)\right\|_1 +\left\|H_{\mathcal{T}^c}\mathbf{x}_e^{\perp}\right\|_1 \leq \epsilon.
\end{equation}
Using the reverse triangle inequality $\|H_{\mathcal{T}}(\mathbf{x}_e^{\perp}-\mathbf{x}_e)\|_1\geq  \big\lvert\|H_{\mathcal{T}}\mathbf{x}_e\|_1-\|H_{\mathcal{T}}\mathbf{x}_e^{\perp}\|_1 \big\rvert \geq \|H_{\mathcal{T}}\mathbf{x}_e\|_1-\|H_{\mathcal{T}}\mathbf{x}_e^{\perp}\|_1$, then \eqref{equ:optimal_projection} implies
\begin{equation}\label{equ:FDIA_proof_2}
    \|H_{\mathcal{T}}\mathbf{x}_e\|_1-\|H_{\mathcal{T}}\mathbf{x}_e^{\perp}\|_1 + \|H_{\mathcal{T}^c}\mathbf{x}_e^{\perp}\|_1 \leq  \epsilon.
\end{equation}
Notice $\sigma_1 = \max_{\mathbf{v}\in \mathbb{R}^n \backslash \{\mathbf{0}\} } \frac{ \|H_{\mathcal{T}}\mathbf{v}\|_1}{\|H_{\mathcal{T}^c}\mathbf{v}\|_1}$, combining \eqref{equ:FDIA_Proof_1} and \eqref{equ:FDIA_proof_2} yields
\begin{equation*} \|H_{\mathcal{T}}\mathbf{x}_e^{\perp}\|_1 - \|H_{\mathcal{T}^c}\mathbf{x}_e^{\perp}\|_1 \geq (\sigma_1-1)\epsilon.
\end{equation*}
Using the topological equivalence between $1$-norm and $2$-norm yields
\begin{equation} \label{equ:first_derivation}
\begin{aligned}
    \sqrt{|\mathcal{T}|}\|H_{\mathcal{T}}\mathbf{x}_e^{\perp}\|_2 -\|H_{\mathcal{T}^c}\mathbf{x}_e^{\perp}\|_2 &\geq (\sigma_1-1)\epsilon, \\
    (\sqrt{|\mathcal{T}|}\overline{\sigma}_{\mathcal{T}} -\underline{\sigma}_{\mathcal{T}^c})\|\mathbf{x}_e^{\perp}\|_2 &\geq (\sigma_1-1)\epsilon, \\
    \|\mathbf{x}_e^{\perp}\|_2 &\geq \frac{(\sigma_1-1)\epsilon}{\sqrt{|\mathcal{T}|}\overline{\sigma}_{\mathcal{T}} -\underline{\sigma}_{\mathcal{T}^c}}. 
\end{aligned}
\end{equation}
Since $|\mathcal{T}|> \frac{\underline{\sigma}_{\mathcal{T}^c}^2}{\overline{\sigma}_{\mathcal{T}}^2}$, the denominator of the right-hand side of \eqref{equ:first_derivation} is positive. Moreover, we can choose a $\mathbf{v}_w \in \mathbb{R}^n$ satisfying $H\mathbf{v}_w = \mathbf{w}$, where $\mathbf{w}$ satisfies the condition in \eqref{equ:violation_NSP}. Therefore there exists an $\mathbf{v}_w \in \mathbb{R}^n$ such that $\sigma_1 \geq \frac{ \|H_{\mathcal{T}}\mathbf{v}_w\|_1}{\|H_{\mathcal{T}^c}\mathbf{v}_w\|_1} >1$. Thus, the upper bound of $\|\mathbf{x}_e^{\perp}\|_2$ in \eqref{equ:first_derivation} is always positive. Furthermore, according to \eqref{equ:decoder} and \eqref{equ:FDIA},
$$ \begin{aligned}
    D(\mathbf{y}_I) &= \argmin_{\mathbf{z}}\|\mathbf{y}_I^{\star} + \mathbf{e}_I - H \mathbf{z}\|_1 \\ & =\argmin_{\mathbf{z}}\|H\mathbf{x}^{\star}+P\begin{bmatrix}H_{\mathcal{T}} \\ 0\end{bmatrix}\mathbf{x}_e-H\mathbf{z}\|_1\\
&=\mathbf{x}^{\star} + \argmin_{\mathbf{z}}\left\|P\begin{bmatrix}H_{\mathcal{T}} \\ 0\end{bmatrix}\mathbf{x}_e-H\mathbf{z}\right\|_1.
\end{aligned}
$$
From \eqref{equ:definition_e_perp}, we have $\mathbf{x}^{\star}-\mathcal{D}(\mathbf{y}_I) = -\mathbf{x}_e^{\perp}$. Then
$$\left\|\mathbf{x}^{\star}-\mathcal{D}(\mathbf{y}_I)\right\|_2 = \left\|\mathbf{x}_e^{\perp}\right\|_2 \geq \frac{(\sigma_1-1)\epsilon}{\sqrt{|\mathcal{T}|}\overline{\sigma}_{\mathcal{T}} -\underline{\sigma}_{\mathcal{T}^c}},
$$
which implies the upper bound of $\alpha$ in \eqref{equ:alpha_L1_FDIA}.

Next, we show that $\|\mathbf{y}_I-H\mathcal{D}(\mathbf{y}_I)\|_2 \leq \epsilon$. Since $\mathbf{x}^{\star}-\mathcal{D}(\mathbf{y}_I) = -\mathbf{x}_e^{\perp}$ and $\mathbf{y}_I = \mathbf{y}_I^{\star} + \mathbf{e}_I$, the residual is

\begin{equation}\label{equ:residual_deduction}
\begin{aligned}
\|\mathbf{y}_I-H\mathcal{D}(\mathbf{y}_I)\|_2 & = 
\|\mathbf{y}_I^{\star} + \mathbf{e}_I-H(\mathbf{x}^{\star}+\mathbf{x}_e^{\perp}))\|_2 \\
& =
\left\|H\mathbf{x}^{\star} + P\begin{bmatrix}H_{\mathcal{T}} \\ 0\end{bmatrix}\mathbf{x}_e - H\mathbf{x}^{\star} - P\begin{bmatrix}H_{\mathcal{T}} \\ H_{\mathcal{T}^c}\end{bmatrix}\mathbf{x}_e^{\perp}\right\|_2 \\
& = \left\|P\left(\begin{bmatrix}H_{\mathcal{T}}\\ H_{\mathcal{T}^c}\end{bmatrix}\mathbf{x}_e^{\perp} - \begin{bmatrix}H_{\mathcal{T}} \\ 0\end{bmatrix}\mathbf{x}_e\right)\right\|_2 \\
&\leq \left\|\begin{bmatrix}H_{\mathcal{T}} \\ H_{\mathcal{T}^c}\end{bmatrix}\mathbf{x}_e^{\perp} - \begin{bmatrix}H_{\mathcal{T}} \\ 0\end{bmatrix}\mathbf{x}_e\right\|_1 \\
    &\leq \|H_{\mathcal{T}}(\mathbf{x}_e^{\perp}-\mathbf{x}_e)\|_1 +\|H_{\mathcal{T}^c}\mathbf{x}_e^{\perp}\|_1.
\end{aligned}
\end{equation}
Combining \eqref{equ:optimal_projection} and \eqref{equ:residual_deduction} yields $\|\mathbf{y}_I-H\mathcal{D}(\mathbf{y}_I)\|_2 \leq \epsilon$.
\end{proof}

\begin{remark}
If $\textsf{null}(H_{\mathcal{T}^c})\backslash \textsf{null}(H_{\mathcal{T}}) \neq \emptyset$, let $\mathbf{v}_n \in \textsf{null}(H_{\mathcal{T}^c})\backslash \textsf{null}(H_{\mathcal{T}})$, then we have $\|H_{\mathcal{T}^c}\mathbf{v}_n\|_1=0$ but $\|H_{\mathcal{T}}\mathbf{v}_n\|_1>0$. Thus, $\sigma_1 \geq \frac{ \|H_{\mathcal{T}}\mathbf{v}\|_1}{\|H_{\mathcal{T}^c}\mathbf{v}\|_1} $ is infinite, which implies that the attack in \eqref{equ:FDIA} is $(\epsilon,\alpha)$-successful for all $\epsilon, \alpha \in \mathbb{R}_{+}$.
\end{remark}

%% file: SICON/revised_sections/5-RESP.tex
In this section, we analyze the resilience of a weighted $\ell_1$ estimator where the weights adaption law includes support prior. With CSP and Row-RIP, different estimation errors are derived. The quantified relationship between the learning prior precision and the estimation error is also clarified and analyzed. 
 
\subsection{Resilient Recovery without Support Prior}
 In this subsection, the moving-horizon $\ell_1$ estimator in \eqref{equ:decoder} is considered. The following theorem gives the conditions for the resilient state recovery by the above decoder.
\begin{theorem}[\textbf{Resilient Recovery with CSP}]\label{Thm0:recovery_Withoutprior}
Consider the measurement model in \eqref{equ:measurement_model}, assume the attacks $\mathbf{e}_i \in \Sigma_k$ for all $i \in I$. If $H \triangleright \textsf{CSP}(Tk)$ with a parameter $\beta \in (0,1)$, the estimation error due to the decoder in \eqref{equ:decoder} can be upper bounded as:
\begin{equation}   \label{equ:error_bound0}
    \|\hat{\mathbf{x}}-\mathbf{x}\|_{2} \leq  \frac{2(1+\beta)}{\underline{\sigma}(1-\beta)}\varepsilon,
 \end{equation} 
 where $\underline{\sigma}$ is the smallest singular value of $H$.
\end{theorem}

\begin{proof}
Since $\hat{\mathbf{x}}$ is the optimal solution of \eqref{equ:decoder},
\begin{equation}\label{equ:optimality_deduction}
    \begin{aligned}
    \|\mathbf{y}-H\hat{\mathbf{x}}\|_1 \leq \|\mathbf{y}-H\mathbf{x}^{\star}\|_1 = \|\mathbf{e}\|_1, \\
    \|\mathbf{y}-H\mathbf{x}^{\star}+H(\mathbf{x}^{\star}-\hat{\mathbf{x}})\|_1 \leq\|\mathbf{e}\|_1, \\
    \|\mathbf{e}+H\tilde{\mathbf{x}}\|_1 \leq \|\mathbf{e}\|_1,
    \end{aligned}
\end{equation}
then decomposing the $1$-norm over $\mathcal{T}$ and $\mathcal{T}^c$ yields 
\begin{equation}
\left\|\mathbf{e}_{\mathcal{T}}+H_{\mathcal{T}}\tilde{\mathbf{x}}\right\|_1 +\left\|\mathbf{e}_{\mathcal{T}^c}+H_{\mathcal{T}^c}\tilde{\mathbf{x}}\right\|_1 \leq \left\|\mathbf{e}_{\mathcal{T}}\|_1+ \|\mathbf{e}_{\mathcal{T}^c}\right\|_1.
\end{equation}
Using the reverse triangle inequality yields 
\begin{equation}
\|\mathbf{e}_{\mathcal{T}}\|_1-\|H_{\mathcal{T}}\tilde{\mathbf{x}}\|_1 -\|\mathbf{e}_{\mathcal{T}^c}\|_1+\|H_{\mathcal{T}^c}\tilde{\mathbf{x}}\|_1 \leq \|\mathbf{e}_{\mathcal{T}}\|_1+ \|\mathbf{e}_{\mathcal{T}^c}\|_1,
\end{equation}
thus, it holds that
\begin{equation*}
\|H_{\mathcal{T}^c}\tilde{\mathbf{x}}\|_1 \leq \|H_{\mathcal{T}}\tilde{\mathbf{x}}\|_1 + 2\|\mathbf{e}_{\mathcal{T}^c}\|_1.
\end{equation*} 
Let $\mathbf{h} = H\tilde{\mathbf{x}}$, based on \eqref{equ:noise_bound}, it follows that
\begin{equation}\label{equ:optimality2}
        \|\mathbf{h}_{\mathcal{T}^c}\|_1 \leq \|\mathbf{h}_{\mathcal{T}}\|_1 + 2\varepsilon,
\end{equation}
Since $H \triangleright \textsf{CSP}(Tk)$, there exists $\beta \in\left(0, 1\right)$ such that $\|\mathbf{h}_{\mathcal{T}}\|_1\le\beta\|\mathbf{h}_{\mathcal{T}^c}\|_1$. Thus,
\begin{equation}\label{equ:SCP_condition1}
    \|\mathbf{h}_{\mathcal{T}}\|_1 \leq \frac{2\beta}{1-\beta}\varepsilon.
\end{equation} 
Then, using the triangle inequality and \eqref{equ:optimality2}, \eqref{equ:SCP_condition1} implies
\begin{equation}\label{equ:inequlity_H}
    \begin{aligned}
\|\mathbf{h}\|_2 &\leq \|\mathbf{h}_{\mathcal{T}}\|_1 +  \|\mathbf{h}_{\mathcal{T}^c}\|_1 \\
&\leq 2\|\mathbf{h}_{\mathcal{T}}\|_1+ 2\varepsilon \\
&\leq \frac{2(1+\beta)}{1-\beta}\varepsilon.
\end{aligned} 
\end{equation}
Finally, combining \eqref{equ:inequlity_H} with $\underline{\sigma}\|\tilde{\mathbf{x}}\|_2 \leq \|\mathbf{h}\|_2$ yields the error bound in \eqref{equ:error_bound1}.
\end{proof}
The following corollary extends the results to a stable error bound based on the Row-RIP condition.
\begin{corollary}[\textbf{Resilient Recovery with Row-RIP}]\label{Thm:recovery_Withoutprior}
Suppose $H \triangleright \textsf{rRIP}((a+1)Tk)$ for some $a>\max{\left(1,\frac{1}{(\underline{\sigma}-1)^2}\right)}$, $a\in \frac{1}{Tk}\mathbb{Z}$. If $k \leq \frac{m}{2a+1}$ and
\begin{equation} \begin{aligned}\label{equ:condition_rip_Withoutprior}
  \delta_{(a+1)Tk} < \frac{a(\underline{\sigma}-\mu_1)^2}{(1+\sqrt{a})^2} -1,
\end{aligned}\end{equation} 
where $\mu_1>0$ and $\underline{\sigma}$ denotes the smallest singular value of $H$, then the estimation error due to the decoder in \eqref{equ:decoder} can be upper bounded as:
\begin{equation}   \label{equ:error_bound1}
    \|\hat{\mathbf{x}}-\mathbf{x}\|_{2} \leq \frac{2\varepsilon}{\mu_1\sqrt{aTk}}.
\end{equation} 
\end{corollary}
\begin{proof}
For the sake of clean presentation, we use both $\mathcal{T}_0$ and $\mathcal{T}$ to represent the actual attack support in the proof. We sort the coefficients of $\mathbf{h}_{\mathcal{T}_0^c}$ in descending order, and let $\mathcal{T}_j, j\in \{1,2,\cdots\}$\footnote{$k\leq \frac{m}{2a+1}$ guarantees at least $\mathcal{T}_j, j\in \{1,2\}$ exist since $\sum_{j=0}^3|\mathcal{T}_j| = (2a+1)Tk\leq Tm$.} denote $j$th support in $\mathbf{h}_{\mathcal{T}_0^c}$ with size $aTk \in \mathbb{Z}$, where $a>1$.
Then,
\begin{equation} 
    \|\mathbf{h}_{\mathcal{T}_{j}}\|_{2} \leq (aTk)^{-1/2}\|\mathbf{h}_{\mathcal{T}_{j}}\|_{\infty}  \leq (aTk)^{-1/2} \|\mathbf{h}_{\mathcal{T}_{j-1}}\|_{1}.
\end{equation} 
Let $\mathcal{T}_{01} = \mathcal{T}_0 \cup \mathcal{T}_1$, we have
\begin{equation} \label{equ:patition_lemma}
\begin{aligned}
    \|\mathbf{h}_{\mathcal{T}_{01}^c}\|_{2} &\leq \sum_{j\geq 2}\|\mathbf{h}_{\mathcal{T}_j}\|_{2}\leq \frac{1}{\sqrt{aTk}}\sum_{j\geq 1}\|\mathbf{h}_{\mathcal{T}_j}\|_{1} \le  \frac{1}{\sqrt{aTk}}\|\mathbf{h}_{\mathcal{T}_0^c}\|_{1}.
    \end{aligned}
\end{equation} 
Following the same proof from \eqref{equ:optimality_deduction} to \eqref{equ:optimality2} and substituting \eqref{equ:patition_lemma} into \eqref{equ:optimality2} yields
\begin{equation} 
\begin{aligned}
    \sqrt{aTk}\|\mathbf{h}_{\mathcal{T}_{01}^c}\|_2 &\leq \|\mathbf{h}_{\mathcal{T}_{0}}\|_1 + 2\varepsilon\\
    & \leq \sqrt{Tk}\|\mathbf{h}_{\mathcal{T}_{0}}\|_2 + 2\varepsilon\\
    & \leq \sqrt{Tk}\|\mathbf{h}_{\mathcal{T}_{01}}\|_2 + 2\varepsilon.
\end{aligned}
\end{equation} 
According to the triangle inequality, 
\begin{equation} 
\begin{aligned}
    \|\mathbf{h}\|_2 \leq \|\mathbf{h}_{\mathcal{T}_{01}^c}\|_2 + \left\|\mathbf{h}_{\mathcal{T}_{01}}\right\|_2 \leq (\frac{1}{\sqrt{a}}+1) \|\mathbf{h}_{\mathcal{T}_{01}}\|_2 + \frac{2\varepsilon}{\sqrt{aTk}}.
\end{aligned}
\end{equation} 
Since $H \triangleright \textsf{rRIP}((a+1)Tk)$ and $\underline{\sigma}\|\tilde{\mathbf{x}}\|_2 \leq \|\mathbf{h}\|_2$, we obtain
\begin{equation} \label{equ:rip_final_bound1}
\begin{aligned}
    \underline{\sigma}\|\tilde{\mathbf{x}}\|_2 \leq \left(\frac{1}{\sqrt{a}}+1\right)\sqrt{1+\delta_{(a+1)Tk}} \|\tilde{\mathbf{x}}\|_2 + \frac{2\varepsilon}{\sqrt{aTk}}.
\end{aligned}
\end{equation} 
Since $\delta_{(a+1)Tk}$ satisfies the condition in \eqref{equ:condition_rip_Withoutprior},
\begin{equation}\label{equ:rip_condition_imply}
    \underline{\sigma} - \left(\frac{1}{\sqrt{a}}+1\right)\sqrt{1+\delta_{(a+1)Tk}} > \mu_1.
\end{equation}
Notice $a>\frac{1}{(\underline{\sigma}-1)^2}$ implies the existence of such $\mu_1>0$. Then \eqref{equ:rip_final_bound1} becomes
\begin{equation}
    \mu_1 \|\tilde{\mathbf{x}}\|_2 \leq \frac{2\varepsilon}{\sqrt{aTk}}
\end{equation}
which implies the estimation error bound in \eqref{equ:error_bound1}.
\end{proof}

\subsection{Resilient Recovery with Support Prior}

\begin{figure}[t!]
    \centering
    \includegraphics[scale=0.55]{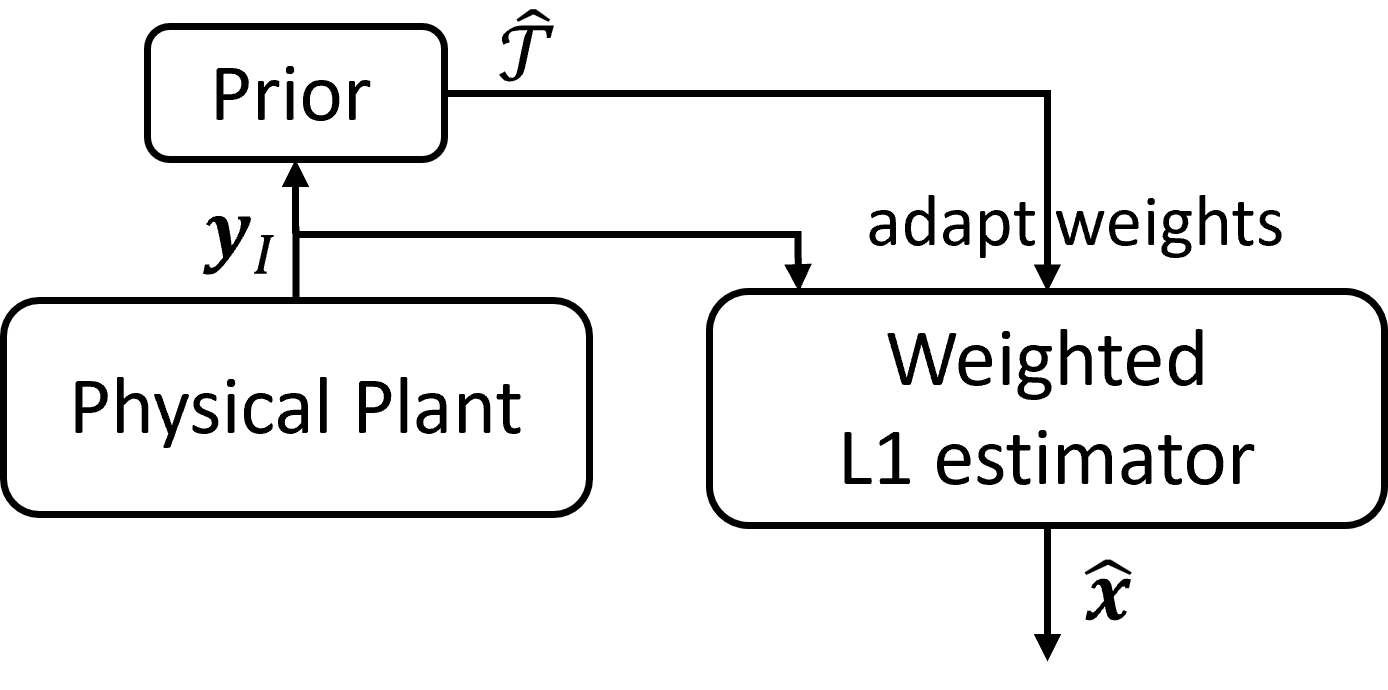}
    \caption{Schematic diagram of resilient estimation scheme with a data-driven prior}
    \label{fig:CLRE}
\end{figure}
In this subsection, we present a resilient estimation scheme with an estimate of attack support, also known as the support prior. A weighted $\ell_1$ estimator is utilized with weight adaption based on support prior. The scheme is illustrated in Figure~\ref{fig:CLRE}. The prior information could be generated from any learning-based attack detector but estimates the attack support
$$\hat{\mathcal{T}} = \{\hat{\mathcal{T}}_i, \cdots, \hat{\mathcal{T}}_{i-T+1}\}.$$
The attack detector performs binary classification on each sensor channel. Let positive class (P) denote the node index that is in the support prior, while negative class (N) denotes the node index that is not in the support prior. Then a receiver operating characteristic (ROC) \cite{green1966signal} table in Table \ref{tab:ROC} could be used to evaluate the prior.
\begin{table}[h!]
    \centering
    \begin{tabular}{|c|c|c|}
      \hline
         $\mathcal{T} \backslash  \hat{\mathcal{T}}$ &  P & N\\
         \hline
         P & TP & FN \\
         \hline
         N & FP & TN \\
         \hline
    \end{tabular}
    \caption{ROC of the prior (TP: true positive, TN: true negative, FP: false positive, FN: false negative)}
    \label{tab:ROC}
\end{table}

\noindent We use a deterministic metric, called the positive predictive value ($\textsf{PPV}$, also called precision), to evaluate the qualification of prior support. The metric is given by
 \begin{equation} \textsf{PPV} \triangleq \frac{TP \text{(true positive)}}{TP+FP\text{(false positive)}} = \frac{|\mathcal{T} \cap  \hat{\mathcal{T}}|}{|\hat{\mathcal{T}}|}.
 \end{equation}
 $\textsf{PPV}$ is chosen as the metric because the positive class is used for the support prior, where minimizing false cases typically enhances the contribution to the resilient estimation scheme. If $\textsf{PPV} > 50\%$, we say the prior is better than the random flip of a fair coin or random guessing. 
 
Next, the corresponding moving-horizon estimator is given as a weighted $\ell_1$ minimization program:
\begin{equation}\label{equ:MHE_weighted_l1_minimization}
    \begin{aligned}
    \textsf{Minimize} &\hspace{2mm} \sum_{j=i-T+1}^i \|\mathbf{y}_j-C\mathbf{x}_j\|_{1,\mathbf{w}(\hat{\mathcal{T}}_j,\omega)}\\
    \textsf{Subject to}& \hspace{2mm} \mathbf{x}_{j+1} -A\mathbf{x}_j =0, \hspace{2mm} j = I-1,
        \end{aligned}
\end{equation}
where $\omega \in\left(0,1\right)$. The optimization problem in \eqref{equ:MHE_weighted_l1_minimization} is equivalent to
 \begin{equation}\label{equ:weighted_l1_minimization}
\begin{aligned}
\Minimize_{\mathbf{z}\in \mathbb{R}^n}\hspace{1mm} \|\mathbf{y}_I - H \mathbf{z}\|_{1, \mathbf{w}(\hat{\mathcal{T}},\omega)},
\end{aligned}    
\end{equation}
where $\mathbf{w}(\hat{\mathcal{T}},\omega) = \begin{bmatrix}\mathbf{w}(\hat{\mathcal{T}_i},\omega) \\ \vdots \\\mathbf{w}(\hat{\mathcal{T}}_{i-T+1},\omega) \end{bmatrix} \in \mathbb{R}^{Tm}$.

Next, a theorem is given to build a quantitative bridge between the prior precision and the estimation error based on the same Row-RIP condition.

\begin{theorem}[\textbf{Resilient Recovery with support prior}]\label{Thm:recovery_Withprior}
Consider the measurement model in \eqref{equ:measurement_model}, assume the attacks $\mathbf{e}_i \in \Sigma_k$ for all $i \in I$ with noise bounded as $\sum_{i\in\mathcal{T}^c}|\mathbf{e}_i|<\varepsilon$. Let $\mathcal{T}$ be the unknown support of attack vector $\mathbf{e}_I$, then assume an estimate $\hat{\mathcal{T}}$ satisfying
\begin{equation} \begin{aligned} \label{equ:support_preformance}
    |\hat{\mathcal{T}}| = \rho |\mathcal{T}| \text{ and }    |\mathcal{T} \cap \hat{\mathcal{T}}| = \textsf{PPV}|\hat{\mathcal{T}}|,
\end{aligned}\end{equation}
where $\rho>0$ and $\textsf{PPV} \in \left(0,\hspace{1mm}1\right)$ is the precision of the estimate. If the following conditions are satisfied:
\begin{enumerate}
    \item $H \triangleright \textsf{rRIP}((a+1)Tk)$ for some $a>\max{\left(1,\frac{1}{(\underline{\sigma}-1)^2}, (1-\textsf{PPV})\rho \right)}, a\in \frac{1}{Tk}\mathbb{Z}$, where $k \leq \frac{m}{2a+1}$,
\item The Row-RIP constant satisfies condition in \eqref{equ:condition_rip_Withoutprior} with same $\mu_1$,
\end{enumerate}
then the estimation error due to the program in \eqref{equ:weighted_l1_minimization} can be upper bounded as:
\begin{equation}\label{equ:error_bound2}
 \|\hat{\mathbf{x}}-\mathbf{x}\|_{2} \leq \frac{2\varepsilon}{\mu_2 \sqrt{aTk}},
\end{equation}
where
\begin{equation}\label{equ:mu2}
    \mu_2 = \mu_1 + \frac{(1-\omega)(1-\sqrt{\kappa})}{\sqrt{a}}\sqrt{1+\delta_{(a+1)Tk}},
\end{equation}
and $\kappa \triangleq 1+\rho- 2\textsf{PPV}\rho$.
\end{theorem}
\begin{proof}
 For a vector $\mathbf{z} \in \mathbb{R}^n$ and two sets $\mathcal{S}_1, \mathcal{S}_2 \subset \{1,2,\dots,n\}$, it holds that
$$ \|\mathbf{z}_{\mathcal{S}_1\cup \mathcal{S}_2}\|_1 = \|\mathbf{z}_{\mathcal{S}_1}\|_1 + \|\mathbf{z}_{\mathcal{S}_2}\|_1 - \|\mathbf{z}_{\mathcal{S}_1\cap \mathcal{S}_2}\|_1.
$$
If $\mathcal{S}_1$ and $\mathcal{S}_2$ are disjoint sets (i.e. $\mathcal{S}_1\cap \mathcal{S}_2 = \emptyset$), then $\|\mathbf{z}_{\mathcal{S}_1\cup \mathcal{S}_2}\|_1 = \|\mathbf{z}_{\mathcal{S}_1}\|_1 + \|\mathbf{z}_{\mathcal{S}_2}\|_1$.
Additionally, we use $\mathcal{T}_0 = \mathcal{T}$ to represent the unknown support of attack vector $\mathbf{e}_I$, thus $|\mathcal{T}_0|\leq Tk$. We also use $\hat{\mathcal{T}}_0 = \hat{\mathcal{T}}$ to represent the estimated support of the attack support. 

Let $\hat{\mathbf{x}}$ be the optimal solution of \eqref{equ:weighted_l1_minimization}, $\tilde{\mathbf{x}} = \mathbf{x} - \hat{\mathbf{x}}$, and $\mathbf{h} = H\tilde{\mathbf{x}}$. Similar to \eqref{equ:optimality_deduction}, we have
\begin{equation}
    \|\mathbf{e} + \mathbf{h}\|_{1,\mathbf{w}(\hat{\mathcal{T}_0},\omega)} \leq \|\mathbf{e}\|_{1,\mathbf{w}(\hat{\mathcal{T}_0},\omega)}.
\end{equation}
By the definition of the weighted $1$-norm, it follows that 
\begin{equation}
\omega \left\|\mathbf{e}_{\hat{\mathcal{T}_0}}+\mathbf{h}_{\hat{\mathcal{T}_0}}\right\|_1 + \left\|\mathbf{e}_{\hat{\mathcal{T}_0}^c}+\mathbf{h}_{\hat{\mathcal{T}_0}^c}\right\|_1 \leq \omega \left\|\mathbf{e}_{\hat{\mathcal{T}_0}}\right\|_1 + \left\|\mathbf{e}_{\hat{\mathcal{T}_0}^c}\right\|_1.
\end{equation}
Then decomposing the $1$-norm over $\mathcal{T}_0$ and $\mathcal{T}_0^c$ on both sides yields 
\begin{equation}
\begin{aligned}
\omega \left\|\mathbf{e}_{\hat{\mathcal{T}_0} \cap \mathcal{T}_0}+\mathbf{h}_{\hat{\mathcal{T}_0}\cap \mathcal{T}_0}\right\|_1 + \omega \left\|\mathbf{e}_{\hat{\mathcal{T}_0} \cap \mathcal{T}_0^c}+\mathbf{h}_{\hat{\mathcal{T}_0}\cap \mathcal{T}_0^c}\right\|_1 + \\
\left\|\mathbf{e}_{\hat{\mathcal{T}_0}^c \cap \mathcal{T}_0}+\mathbf{h}_{\hat{\mathcal{T}_0}^c \cap \mathcal{T}_0}\right\|_1 + \left\|\mathbf{e}_{\hat{\mathcal{T}_0}^c \cap \mathcal{T}_0^c}+\mathbf{h}_{\hat{\mathcal{T}_0}^c \cap \mathcal{T}_0^c}\right\|_1 \leq \\
\omega \left\|\mathbf{e}_{\hat{\mathcal{T}_0}\cap \mathcal{T}_0}\right\|_1 + \left\|\mathbf{e}_{\hat{\mathcal{T}_0}\cap \mathcal{T}_0^c}\right\|_1 + \left\|\mathbf{e}_{\hat{\mathcal{T}_0}^c\cap \mathcal{T}_0}\right\|_1 + \left\|\mathbf{e}_{\hat{\mathcal{T}_0}^c\cap \mathcal{T}_0^c}\right\|_1.
\end{aligned}
\end{equation}
By using the reverse triangle inequality and arranging the inequality, we obtain
\begin{equation}\label{equ:54_proof1}
    \begin{aligned}\omega \left\|\mathbf{h}_{\hat{\mathcal{T}_0}\cap \mathcal{T}_0^c}\right\|_1 + \left\|\mathbf{h}_{\hat{\mathcal{T}_0}^c\cap \mathcal{T}_0^c}\right\|_1 \leq \\ \left\|\mathbf{h}_{\hat{\mathcal{T}_0}^c\cap \mathcal{T}_0}\right\|_1  +\omega\left\|\mathbf{h}_{\hat{\mathcal{T}_0}\cap \mathcal{T}_0}\right\|_1
+2(\|\mathbf{e}_{\hat{\mathcal{T}_0}^c\cap \mathcal{T}_0^c}\left\|_1+\|\mathbf{e}_{\hat{\mathcal{T}_0}\cap \mathcal{T}_0^c}\right\|_1).
\end{aligned}
\end{equation}
Next, adding and subtracting $\omega\|\mathbf{h}_{\hat{\mathcal{T}_0}^c\cap \mathcal{T}_0^c}\|_1$ on the left side of \eqref{equ:54_proof1}, and $\omega \|\mathbf{h}_{\hat{\mathcal{T}_0}^c\cap \mathcal{T}_0}\|_1$ and $\omega\|\mathbf{e}_{\hat{\mathcal{T}_0}^c\cap \mathcal{T}_0^c}\|_1$ on the right yields 
\begin{equation}\label{equ:54_proof2}
\begin{aligned}
\omega \left\|\mathbf{h}_{\mathcal{T}_0^c}\right\|_1 +(1-\omega)\left\|\mathbf{h}_{\hat{\mathcal{T}_0}^c\cap \mathcal{T}_0^c}\right\|_1  \leq \\  (1-\omega) \left\|\mathbf{h}_{\hat{\mathcal{T}_0}^c\cap \mathcal{T}_0}\right\|_1 +\omega \left\|\mathbf{h}_{\mathcal{T}_0}\right\|_1
+2(\omega \|\mathbf{e}_{\mathcal{T}_0^c}\|_1 +(1-\omega) \left\|\mathbf{e}_{\hat{\mathcal{T}_0}^c\cap \mathcal{T}_0^c}\right\|_1 ).
\end{aligned}
\end{equation}
 Again, adding and subtracting $(1-\omega)\|\mathbf{h}_{\hat{\mathcal{T}_0}\cap \mathcal{T}_0^c}\|_1$ on the left side of \eqref{equ:54_proof2} and substituting $\sum\limits_{i\in\mathcal{T}_0^c}|\mathbf{e}_i|<\varepsilon$ yields
\begin{equation}\label{equ:final_before1}
\begin{aligned}
        \|\mathbf{h}_{\mathcal{T}_0^c}\|_1 \leq \omega  \|\mathbf{h}_{\mathcal{T}_0}\|_1 +(1-\omega)(\left\|\mathbf{h}_{\hat{\mathcal{T}_0}^c\cap \mathcal{T}_0}\right\|_1+ \left\|\mathbf{h}_{\hat{\mathcal{T}_0}\cap \mathcal{T}_0^c}\right\|_1) 
        +2\varepsilon.  
        \end{aligned}
\end{equation}
Since $\hat{\mathcal{T}_0}^c \cap \mathcal{T}_0$ and $\hat{\mathcal{T}_0} \cap \mathcal{T}_0^c$ are disjoint and $(\hat{\mathcal{T}_0}^c \cap \mathcal{T}_0) \cup (\hat{\mathcal{T}_0} \cap \mathcal{T}_0^c) = \mathcal{T}_0 \cup \hat{\mathcal{T}}_0 \setminus \hat{\mathcal{T}_0} \cap \mathcal{T}_0$, we rewrite \eqref{equ:final_before1} as
\begin{equation}\label{equ:final_before}
\begin{aligned}
        \|\mathbf{h}_{\mathcal{T}_0^c}\|_1 \leq \omega  \|\mathbf{h}_{\mathcal{T}_0}\|_1 +(1-\omega)\left\|\mathbf{h}_{ \mathcal{T}_0 \cup \hat{\mathcal{T}}_0 \setminus \hat{\mathcal{T}}_0 \cap \mathcal{T}_0}\right\|_1 
        +2\varepsilon.  
        \end{aligned}
\end{equation}

 We sort the coefficients of $\mathbf{h}_{\mathcal{T}_0^c}$ in descending order, and let $\mathcal{T}_j, j\in \{1,2,\cdots\}$ denote $j$th support in $\mathbf{h}_{\mathcal{T}_0^c}$ with size $aTk \in \mathbb{Z}$, where $a>1$. Since $\mathbf{h}_{\mathcal{T}_1}$ corresponds to the biggest coefficients in $\mathbf{h}_{\mathcal{T}_0^c}$, and $|\hat{\mathcal{T}}_0 \setminus \hat{\mathcal{T}}_0 \cap \mathcal{T}_0| = (1-\textsf{PPV})\rho Tk \leq aTk = |\mathcal{T}_1|$, we have $\left\|\mathbf{h}_{\hat{\mathcal{T}}_0 \setminus \hat{\mathcal{T}}_0 \cap \mathcal{T}_0}\right\|_2 \leq \left\|\mathbf{h}_{\mathcal{T}_1}\right\|_2$, then $\left\|\mathbf{h}_{ \mathcal{T}_0 \cup \hat{\mathcal{T}}_0 \setminus \hat{\mathcal{T}}_0 \cap \mathcal{T}_0}\right\|_2 \leq \left\|\mathbf{h}_{\mathcal{T}_{01}}\right\|_2$. Additionally, according to the topological equivalence of vector norm in \eqref{equ:topo_equi_norm}, since $| \mathcal{T}_0 \cup \hat{\mathcal{T}}_0 \setminus \hat{\mathcal{T}}_0 \cap \mathcal{T}_0| = \kappa Tk$, it follows that
 \begin{equation}\label{equ:connection_proof}
     \left\|\mathbf{h}_{ \mathcal{T}_0 \cup \hat{\mathcal{T}}_0 \setminus \hat{\mathcal{T}}_0 \cap \mathcal{T}_0}\right\|_1 \leq \sqrt{\kappa Tk}\left\|\mathbf{h}_{ \mathcal{T}_0 \cup \hat{\mathcal{T}}_0 \setminus \hat{\mathcal{T}}_0 \cap \mathcal{T}_0}\right\|_2 \leq \sqrt{\kappa Tk} \left\|\mathbf{h}_{\mathcal{T}_{01}}\right\|_2.
 \end{equation}
Similarly, 
\begin{equation}\label{equ:connection_proof2}
    \|\mathbf{h}_{\mathcal{T}_0}\|_1 \leq \sqrt{Tk}\|\mathbf{h}_{\mathcal{T}_0}\|_2 \leq \sqrt{Tk}\|\mathbf{h}_{\mathcal{T}_{01}}\|_2.
\end{equation}
Substituting \eqref{equ:connection_proof2} and \eqref{equ:connection_proof} into \eqref{equ:final_before} yields
\begin{equation}\label{equ:result_optimization_weighted}
\begin{aligned}
        \left\|\mathbf{h}_{\mathcal{T}_0^c}\right\|_1 \leq \left[\omega \sqrt{Tk} +(1-\omega)\sqrt{\kappa Tk}\right]\left\|\mathbf{h}_{ \mathcal{T}_{01}}\right\|_2 
        +2\varepsilon.  
        \end{aligned}
\end{equation}
Substituting \eqref{equ:patition_lemma} into \eqref{equ:result_optimization_weighted} yields
\begin{equation}
\begin{aligned}
        \sqrt{aTk} \|\mathbf{h}_{ \mathcal{T}_{01}^c}\|_2 \leq \left[\omega \sqrt{Tk} +(1-\omega)\sqrt{\kappa Tk}\right]\|\mathbf{h}_{ \mathcal{T}_{01}}\|_2 
        +2\varepsilon.
        \end{aligned}
\end{equation}
 By applying the triangle inequality, we obtain
\begin{equation}
\begin{aligned}
        \|\mathbf{h}\|_2 \leq \left(\frac{\omega+(1-\omega)\sqrt{\kappa}}{\sqrt{a}}+1\right) \|\mathbf{h}_{ \mathcal{T}_{01}}\|_2
        +\frac{2\varepsilon}{\sqrt{aTk}}.
        \end{aligned}
\end{equation}
Since $H \triangleright \textsf{rRIP}((a+1)Tk)$ and $\underline{\sigma}\|\tilde{\mathbf{x}}\|_2 \leq \|\mathbf{h}\|_2$, we obtain
\begin{equation}
        \underline{\sigma}\|\tilde{\mathbf{x}}\|_2 \leq \left(\frac{\omega+(1-\omega)\sqrt{\kappa}}{\sqrt{a}}+1\right)\sqrt{1+\delta_{(a+1)Tk}} \|\tilde{\mathbf{x}}\|_2 
        + \frac{2\varepsilon}{\sqrt{aTk}}.
\end{equation}
Subtracting $\left(\frac{1}{\sqrt{a}}+1\right)\sqrt{1+\delta_{(a+1)Tk}}\|\tilde{\mathbf{x}}\|_2$ on both sides yields
$$ \begin{aligned}
    \left[\underline{\sigma} - \left(\frac{1}{\sqrt{a}}+1\right)\sqrt{1+\delta_{(a+1)Tk}} \right]\|\tilde{\mathbf{x}}\|_2 \leq \frac{\omega+(1-\omega)\sqrt{\kappa}-1}{\sqrt{a}}\sqrt{1+\delta_{(a+1)Tk}} \|\tilde{\mathbf{x}}\|_2 \\
    + \frac{2\varepsilon}{\sqrt{aTk}}.
\end{aligned}
$$
Since $\delta_{(a+1)Tk}$ satisfies the condition in \eqref{equ:condition_rip_Withoutprior}, the inequality in \eqref{equ:rip_condition_imply} holds. Then
$$\begin{aligned}
    \mu_1 \|\tilde{\mathbf{x}}\|_2 \leq \frac{\omega+(1-\omega)\sqrt{\kappa}-1}{\sqrt{a}}\sqrt{1+\delta_{(a+1)Tk}} \|\tilde{\mathbf{x}}\|_2 
    + \frac{2\varepsilon}{\sqrt{aTk}}.
\end{aligned} 
$$
Let $\mu_2$ be given as in \eqref{equ:mu2}, then 
$$ \mu_2 \|\tilde{\mathbf{x}}\|_2 \leq \frac{2\varepsilon}{\sqrt{aTk}} 
$$
which implies the estimation error bound in \eqref{equ:error_bound2}.
\end{proof}

\begin{remark}\label{rmk:compare_error_bounds}
By comparing the estimation error bound in \eqref{equ:error_bound1} and \eqref{equ:error_bound2}, resilient recovery has a smaller estimation error by adding the measurement prior if $\mu_2>\mu_1$, which is guaranteed by $ \kappa =\rho+1-2\textsf{PPV}\rho<1$. This condition is equivalent to $ \textsf{PPV} > 50\%$. In other words, if the prior is at least better than the random flip of a fair coin, then the resilient recovery with prior outperforms the resilient recovery without prior.
\end{remark}

\begin{figure}[h!]
    \centering
    \includegraphics[scale=0.5]{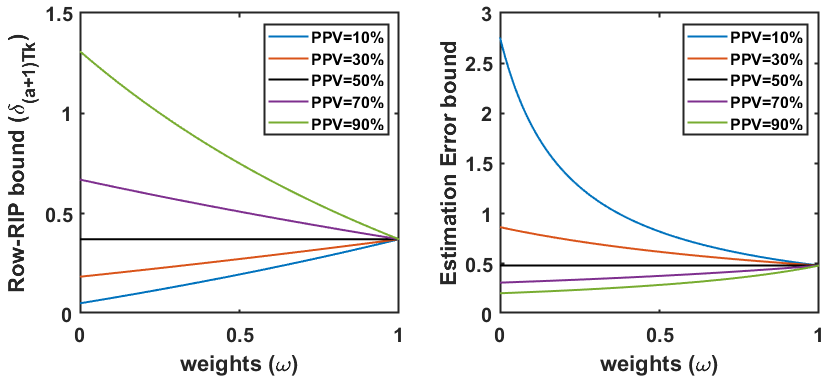}
    \caption{The effects of weight setting $\omega$ of the weighted $\ell_1$ decoder in \eqref{equ:weighted_l1_minimization} and the precision $\textsf{PPV}$ of the prior on the performance of resilient state recovery. (Other variables are set as $\rho=1$, $\underline{\sigma}=2$, $a=2$, $Tk = 50$)}
    \label{fig:weight-analysis}
\end{figure}

To visualize the effect of the weight $\omega$ of the weighted $\ell_1$ decoder in \eqref{equ:weighted_l1_minimization} and the precision $\textsf{PPV}$ of the prior on the performance of the state recovery, a plot of the estimation error bound in \eqref{equ:error_bound2} and the Row-RIP constant's upper bound in \eqref{equ:condition_rip_Withoutprior} is shown in Fig.~\ref{fig:weight-analysis}. $\textsf{PPV} =50\%$ is a boundary where the weight setting does not improve or reduce the performance of the state recovery. Higher precision of the prior leads to better performance of the state recovery (less strict Row-RIP condition and smaller estimation error). Moreover, if the prior is better than the random flip of a fair coin, smaller weights correspond to better resilient recovery performance. If the prior is worse than the random flip of a fair coin, it is best to choose $\omega=1$ which reduces the weighted $\ell_1$ estimation to $\ell_1$ estimation.

%% file: SICON/revised_sections/6-simulation.tex
\begin{figure}[h!]
    \centering
    \includegraphics[scale=0.5]{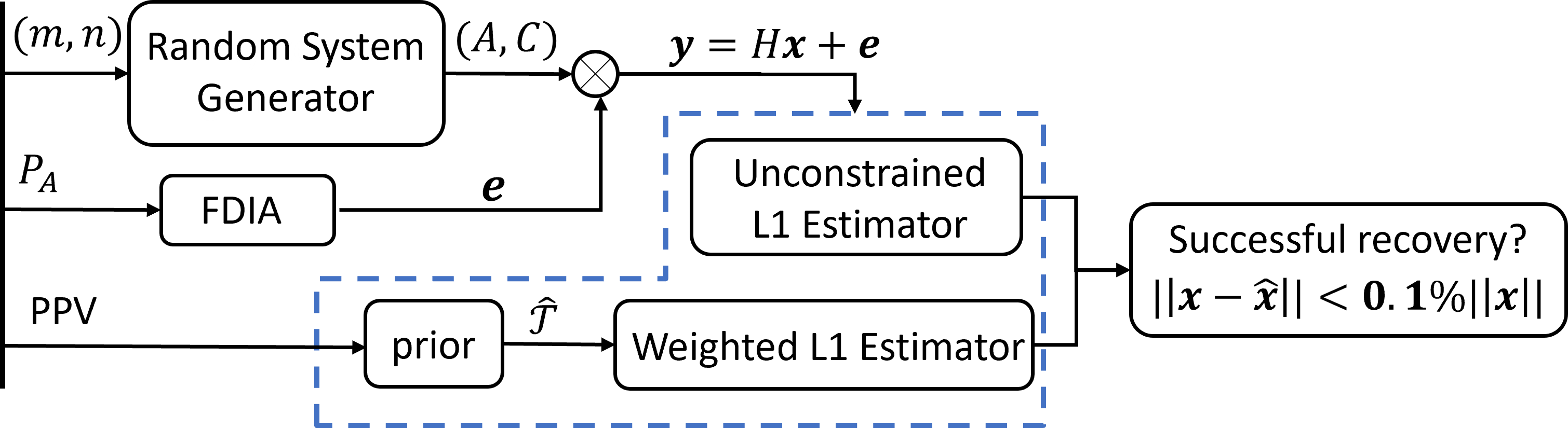}
    \caption{Block diagram depiction of the Monte-Carlos test scenario. ($(m,n)$is the system dimension pair, $m$ is the dimension of measurements, $n$ is the dimension of states, $P_A$ is the attack-to-normal-signal ratio, $(A,C)$ is system dynamic pair, $\textsf{PPV}$ is the precision of the prior)}
    \label{fig:sim_scheme}
\end{figure}

In this section, we present Monte-Carlo tests for the developed resilient state recovery schemes: weighted $\ell_1$ estimator with prior, and $\ell_1$ estimator without prior, as shown in Fig.~\ref{fig:sim_scheme}. The resilient estimation scheme is tested on random linear systems with different attack-to-normal-signal ratios and different prior precision. 

The simulations are processed using the following procedure: 
\begin{enumerate}
    \item Select system dimension pair $(m,n)$, where $m$ is the size of measurements, $n$ is the size of states, and $m > n$. Then generate an observable LTI model $\mathbf{x}_{i+1} = A\mathbf{x}_i, \hspace{2mm} \mathbf{y}_i = C\mathbf{x}_i$, represented by pair $(A,C)$, where $A \in {\mathbb R}^{n \times n}$ and $C \in {\mathbb R}^{m \times n}$ are sampled with independent Gaussian entries \cite{candes2005decoding}. Although the RIP conditions in the previous sections are NP-hard to verify, it has been shown that Gaussian matrices with independent identically distributed entries satisfy the RIP condition with overwhelming probabilities \cite{candes2005decoding, donoho2006most}.
    \item Select an attack-to-normal-signal ratio $P_A$, and find attack location indexed in $\mathcal{T}$ randomly, where $|\mathcal{T}|\leq mP_A$. Then generate attack vectors by \eqref{equ:FDIA_optimization}, and inject it into the selected attack location.
    \item For the unconstrained $\ell_1$ estimator, use linear programming to solve the $\ell_1$ minimization problem \eqref{equ:decoder}, such that
     \begin{equation} \min_{t,\mathbf{x}}\hspace{1mm} \mathbf{1}^{\top} t, \hspace{2mm} \text{subject to}: -t \leq \mathbf{y}_I - H\mathbf{x} \leq t.
     \end{equation}
    \item For the weighted $\ell_1$ estimator, simulate the prior information with a particular $\textsf{PPV}$ by using a Bernoulli distributed agreement model:
     \begin{equation} \label{equ:uncertainty_model}
     \mathbf{q}_{i}=\epsilon_{i}\hat{\mathbf{q}}_{i}+(1-\epsilon_{i})(1-\hat{\mathbf{q}}_{i}),
 \end{equation}
 where $\epsilon_i \sim \mathcal{B}(1,\mathbf{p}_i)$ models the agreement between the estimate support prior $\hat{\mathcal{T}}$ and the actual attack support  $\mathcal{T}$, its probability $\mathbf{p}_i = E[\epsilon_i] = \textsf{Pr}\left\{\epsilon_i = 1\right\}$ approximates the given precision of the support prior $\textsf{PPV}$. Additionally, $\mathbf{q} \in \{0, 1\}^{n}$ is an indicator vector of $\mathcal{T} \subseteq \{1,\dots, n\}$ defined element-wise as:
 \begin{equation} \label{equ:support}
     \mathbf{q}_{i}=\left\{
     \begin{array}{lr}
     0 &\text{if}\hspace{0.2cm} i \in \mathcal{T} \\
     1 &\text{otherwise}.
     \end{array}
     \right.
 \end{equation}   
 Then the weighted $\ell_1$ minimization problem in \eqref{equ:MHE_weighted_l1_minimization} is solved by
    \begin{equation}\label{equ:linprog_weighted_L1}
        \min_{t,\mathbf{x}}\hspace{1mm} \mathbf{w}(\hat{\mathcal{T}},\omega)^{\top} t, \hspace{2mm} \text{subject to}: -t \leq \mathbf{y}_I - H\mathbf{x} \leq t,
    \end{equation} 
    where $\mathbf{w}(\hat{\mathcal{T}},\omega)$ is the weight vector defined as \eqref{equ:weight_definition} with the support prior $\hat{\mathcal{T}}^c$. According to Fig.~\ref{fig:weight-analysis}, we set $\omega=0.01$ when the precision is above $50\%$, and $\omega = 0.99$ when the precision is lower than $50\%$.
    \item Establish a performance evaluation criteria: we determined successful state recovery if the absolute estimation error is less than $0.1\%$ of the nominal absolute value of states. To calculate the success percentage, $10000$ trials of simulation are performed for different resilient estimation schemes.
    \item Repeat the described simulation procedure for various settings of the system dimension $(m,n)$, the attack-to-normal signal ratio $P_A$, and the precision $\textsf{PPV}$.
\end{enumerate}  

\begin{figure}[h!]
    \centering
    \includegraphics[scale=0.48]{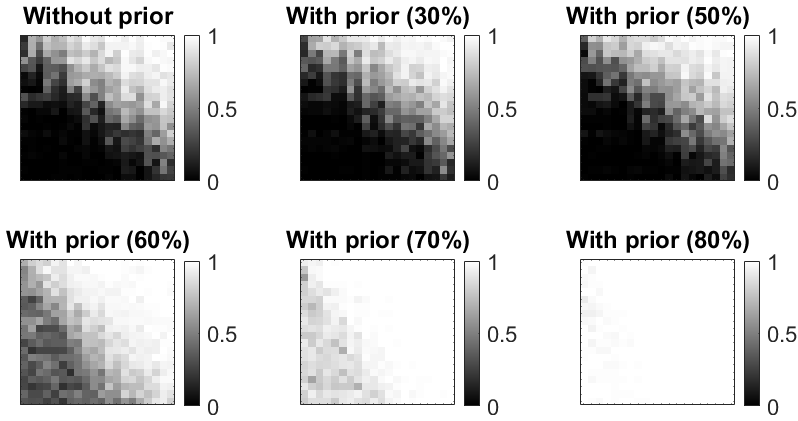}
    \caption{The success ratio of resilient recovery (from black $0\%$ to white $100\%$) of unconstrained $\ell_1$ estimator (without prior) and weighed $\ell_1$ estimator with the prior (with prior) on $400$ different dimensional LTI systems (x-axis: measurement dimension $m_{\text{new}} \in [10,86]$, y-axis: state dimension $n \in [5,43]$) with the attack-to-normal-signal ratio $P_A=40\%$. Different precision, from $30\%$ to $80\%$, of the prior are simulated.}
    \label{fig:result1}
\end{figure}

\begin{figure}[h!]
    \centering
    \includegraphics[scale=0.48]{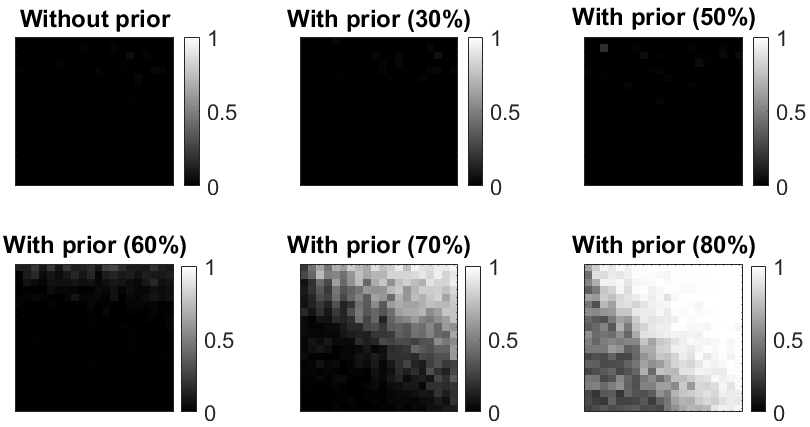}
    \caption{The success ratio of resilient recovery (from black $0\%$ to white $100\%$) of unconstrained $\ell_1$ estimator (without prior) and weighed $\ell_1$ estimator with the prior (with prior) on $400$ different dimensional LTI systems (x-axis: measurement dimension $m_{\text{new}}\in[10,86]$, y-axis: state dimension $n \in [5,43]$) with the attack-to-normal-signal ratio $P_A=60\%$. From left to right, the mean value of agreement probability (positively correlated to the precision of the prior) is $60\%$, $70\%$, and $80\%$ respectively.}
    \label{fig:result2}
\end{figure}

The results in Fig.~\ref{fig:result1} and Fig.~\ref{fig:result2} show the impact of system dimension $(m,n)$ (especially the relationship between $m$ and $n$), the prior precision $\textsf{PPV}$ and the attack-to-normal-signal ratio $P_A$ on the performance of resilient recovery using the unconstrained $\ell_1$ estimator and the weighted $\ell_1$ estimator. The x-axis represents the size of measurements $m \in [10,86]$, while the y-axis represents the dimension of states $n \in [5,43]$. Each small rectangular block in the subplots represents one of $400$ different systems, each with a different dimension $(m_{\text{new}},n)$, where $m_{\text{new}} = \max{(m,n)}$\footnote{$m_{\text{new}}\geq n$ must be guaranteed otherwise state recovery could not be accomplished.}. Thus, the lower left corner of each subplot corresponds to the smallest measurement redundancy, while the upper right corner represents the largest measurement redundancy. We define successful recovery as the case where the absolute estimation error is less than $0.1\%$ of the magnitude of the nominal state. By running $10000$ trials of simulation for each system with different dimensions, the grey level in the plots indicates the ratio of the number of successful estimations to $10000$. The black color indicates $0\%$ successful cases, while the white color indicates $100\%$ successful cases. 

In Fig.~\ref{fig:result1}, $40\%$ measurements are attacked and the precision of the prior is simulated from $30\%$ to $80\%$. The grey level becomes lighter from the lower left corner to the upper right corner in all subplots. This change indicates that larger measurement redundancy enables better resiliency of the system. There is an approximate boundary between the dark grey portion and the light grey portion. For the unconstrained $\ell_1$ estimator, the boundary is around $m=2n$, which is consistent with literature \cite{fawzi2014secure, shoukry2015event}. As shown in the first row, when the precision of the prior is less than or equal to $50\%$, the estimation performance of the weighted $\ell_1$ estimator (with prior) is similar to the unconstrained $\ell_1$ estimator (without prior). As shown in the second row, the boundary moves to the lower left corner more when we use the weighted $\ell_1$ estimation with the prior of higher precision. This shift indicates that if the prior is better than the random flip of a fair coin, less measurement redundancy is required for successful state recovery. By comparing all weighted $\ell_1$ estimators with the prior of different precision, a positive correlation between the precision of the prior and the resilient estimation performance is observed as discussed in Theorem~\ref{Thm:recovery_Withprior}. 

In Fig.~\ref{fig:result2}, we increase the attck-to-normal signal ratio to $60\%$. As proved in literature \cite{fawzi2014secure,shoukry2015event} and in Theorem~\ref{Thm:recovery_Withoutprior}, the 2s-observability assumption is violated so that the unconstrained $\ell_1$ estimator fails. When the prior's precision is less than $50\%$, the estimation performance cannot be improved by including the support prior $\hat{\mathcal{T}}$. Thus, the subplots in the first row are all black. But the second row immediately shows the improvement of the estimation performance by including the support prior $\hat{\mathcal{T}}$ with higher precision than $50\%$.

\begin{figure}[h!]
    \centering
    \includegraphics[scale=0.48]{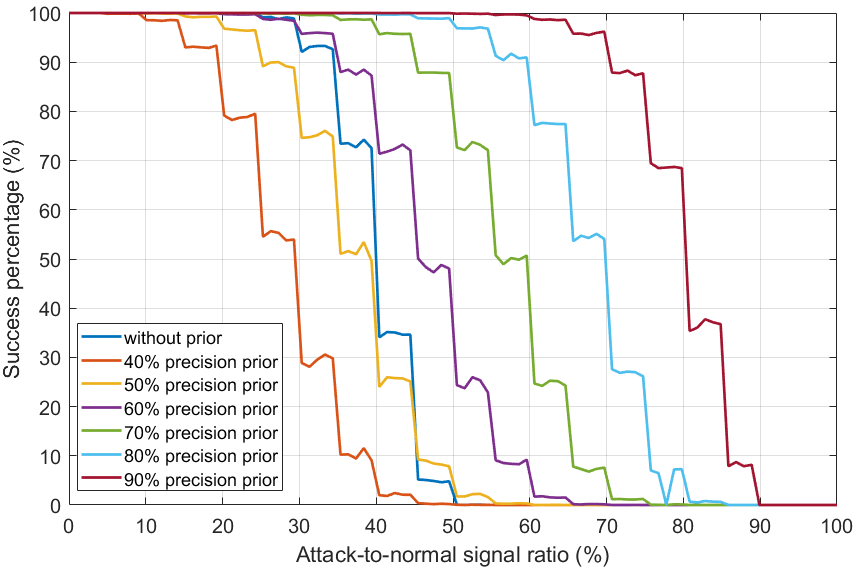}
    \caption{A comparison of estimation performance with different attack-to-normal-signal ratio $P_A \in [0,1]$ between unconstrained $\ell_1$ estimator (without prior) and weighted $\ell_1$ estimator with the prior (with different precision of prior). The precision of the prior is ranging from $40\%$ to $90\%$, and the system dimension is $m=2n, n=10$.}
    \label{fig:result3}
\end{figure}

Furthermore, we plot the success ratio of the resilient recovery of both the unconstrained $\ell_1$ estimator and the weighted $\ell_1$ estimator with the prior for recovering states from attacked measurements, as a function of the attack-to-normal-signal ratio in Fig~\ref{fig:result3}. The S-shape curve between the success ratio of state recovery and the attack percentage indicates that the resilient recovery will maintain a high success ratio when a particular percentage of sensors remain safe. Beyond that point, the success ratio of state recovery rapidly declines. Flat portions of the curves in Fig~\ref{fig:result3} indicate that the success ratio of state recovery remains largely unchanged when there are only minor differences in the number of attacked channels.

Conversely, the dark blue line, representing the unconstrained $\ell_1$ estimator without prior, shows the exact limitation proven in \cite{fawzi2014secure, shoukry2015event}: more than $50\%$ sensors are required to be safe. However, this limitation is relaxed by using a weighted $\ell_1$ estimator with prior as long as the precision of the prior is greater than $50\%$, as proved in Theorem~\ref{Thm:recovery_Withprior}. If the precision of prior is less than $50\%$ (orange line) or lies on the boundary (yellow line), the resilient recovery performance may not be better than the unconstrained $\ell_1$ estimator without prior. This Monte-Carlo simulation confirms the proven limitation of the unconstrained $\ell_1$ estimator in literature and validates the quantitative correlation between the precision of the prior and the resilient recovery performance outlined in this paper.

%% file: SICON/revised_sections/7-case_study.tex
In this section, we use  a linear control system on the IEEE 14-bus system\footnote{\url{https://labs.ece.uw.edu/pstca/pf14/pg_tca14bus.htm}} to validate our theoretical claims. The bus system has $14$ buses and $5$ generators, and each bus in the network is assumed to be equipped with IIoT measurement devices with the capability of measuring active power injections and flow measurements. The detailed physical system construction can be found in \cite{anubi2019enhanced}, which includes the plant model, proportional and integral (PI) regulator, and measurement model. 
The goal of the resilient observer is to estimate $10$ states, including $5$ generator rotor angles and $5$ generator frequencies, correctly from $19$ measurements including $5$ generator frequencies and net powers at $14$ buses.

\begin{figure}[!h]
\begin{center}
\includegraphics[scale=0.5]{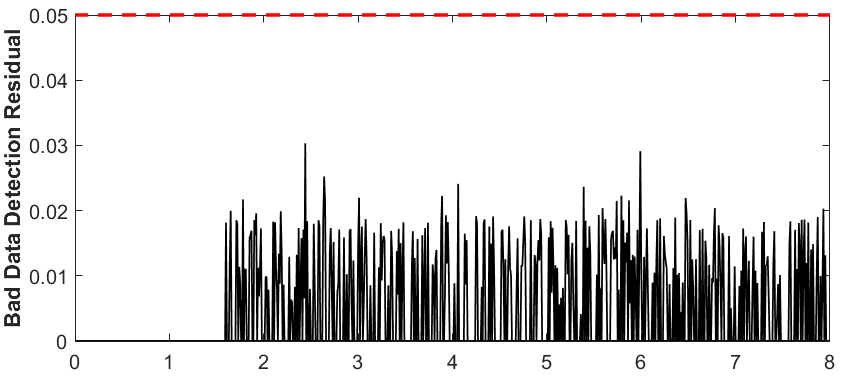}
\end{center}
\caption{Bad data detection result (the residual threshold is set as $0.05$)}
\label{BBD} 
\end{figure}

\begin{figure*}[!h]
\begin{center}
\includegraphics[width=1\textwidth]{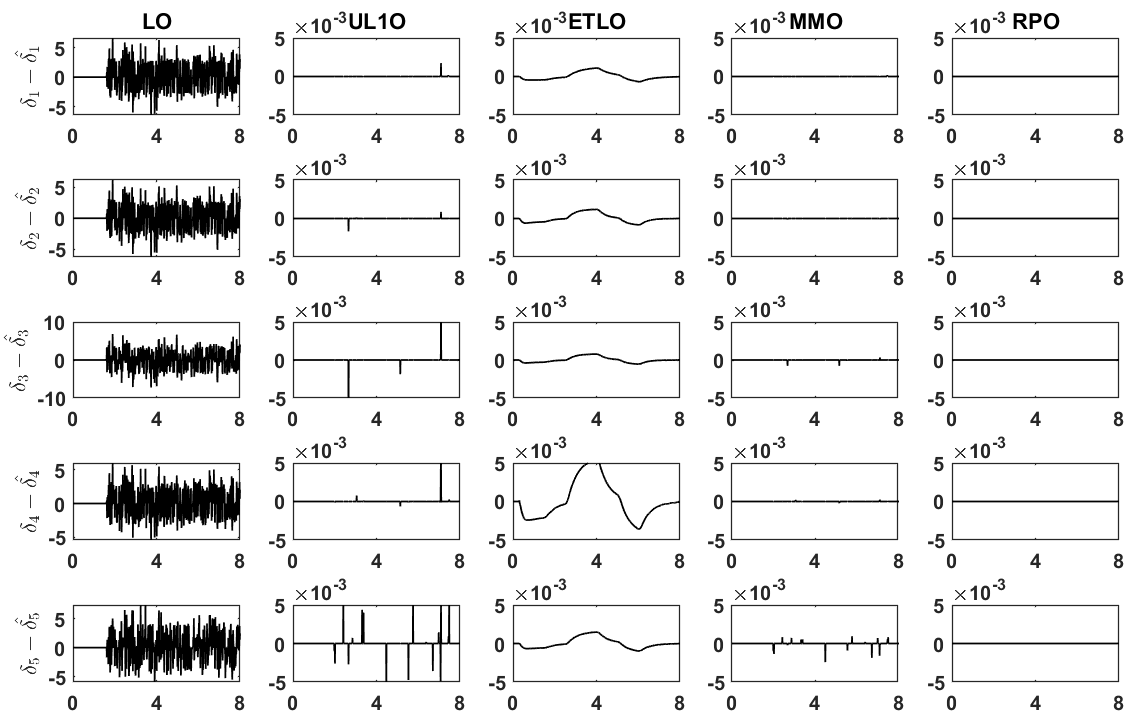}
\end{center}
\caption{A comparison of the estimation errors of generator rotor angles by 5 observers (LO: Luenberger observer, UL1O: unconstrained $\ell_1$ observer, MMO: multi-model observer, ETLO: event-triggered Luenberger observer, RPO: the proposed resilient pruning observer).}
\label{Comparison_all} 
\end{figure*}

We implement the attacks designed in \eqref{equ:FDIA}, $20 \%$ of the measurements are compromised, so $|\mathcal{T}| = 19 \times 0.2 \approx 4$ sensors are randomly selected to be attacked. Fig.~\ref{BBD} indicates the attack designed in \eqref{equ:FDIA} is capable of bypassing the residual-based bad data detector. As shown in Fig.~\ref{Comparison_all}, an observer without resilient design, such as a Luenberger observer, is unstable under such attacks. Although such attacks are too strong to be implemented in practice due to the assumption of knowing the system's model, it is sufficient for the comparison of the resiliency of resilient observers in this application example since resilient estimators are mainly affected by the number of attacked sensors instead of the design of attacks.

\begin{figure}
    \centering
    \includegraphics[width=0.8\linewidth]{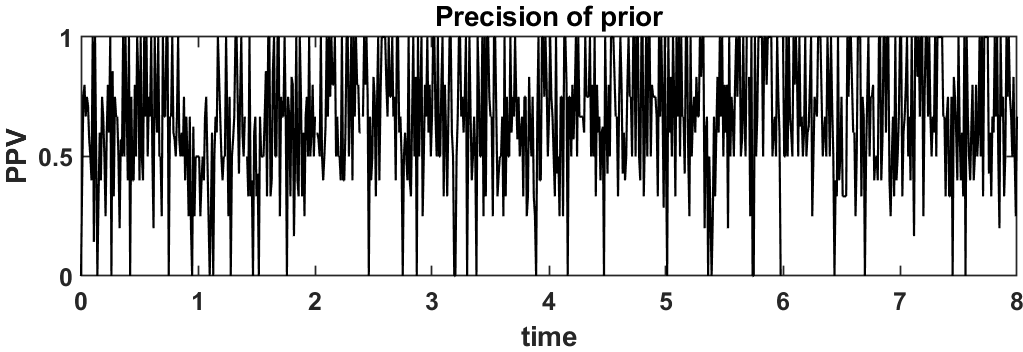}
    \caption{The precision of support prior generated by the GPR-based attack localizer for the power grid (The mean of precision is $65.5\%$)}
    \label{fig:prior_precision}
\end{figure}

\begin{figure}[!t]
\begin{center}
\includegraphics[scale=0.4]{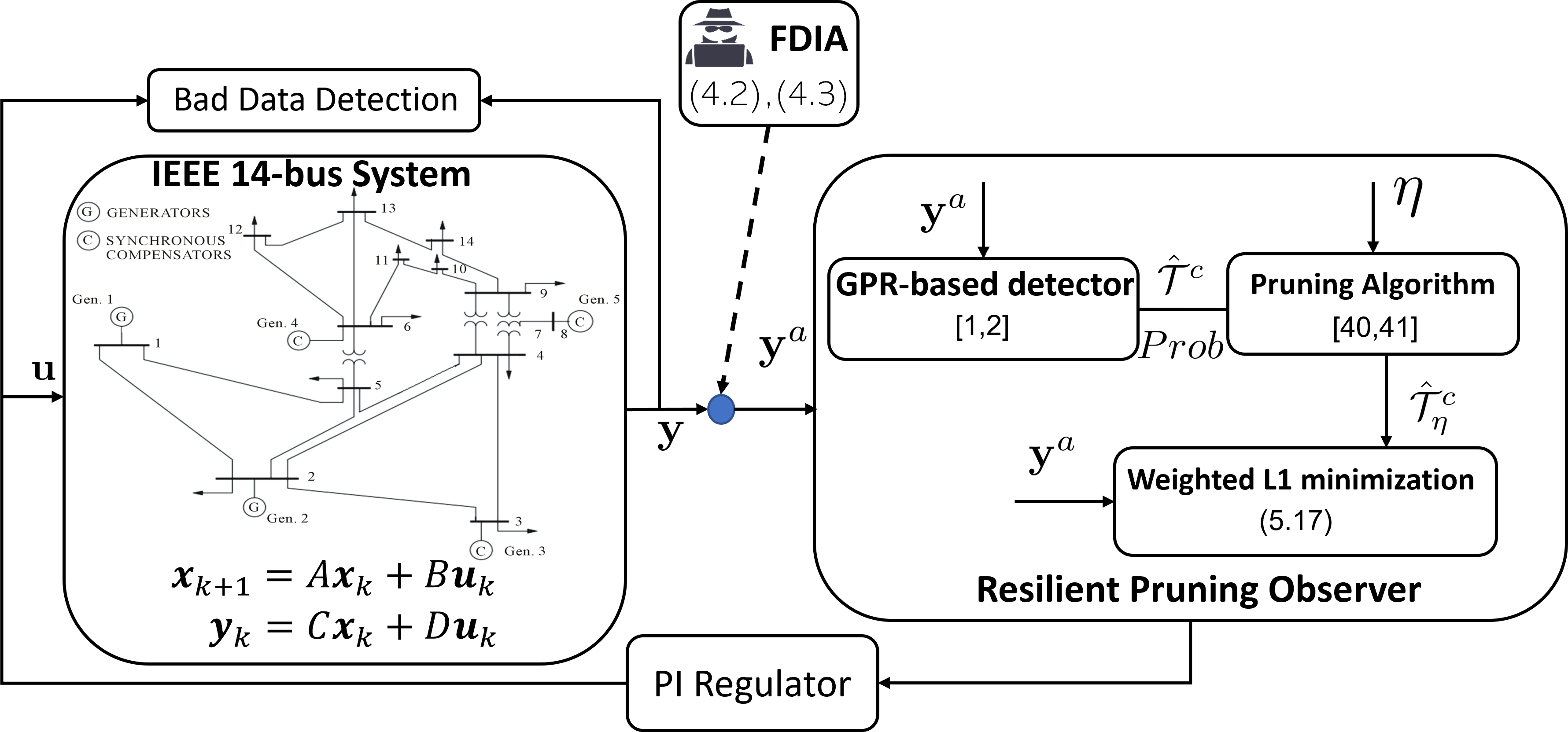}
\end{center}
\caption{Block diagram depiction of the implementation of the resilient pruning observer (RPO) on IEEE 14 bus system (FDIA: false data injection attack).}
\label{Resilient_pruning_simulation} 
\end{figure}

We used a Gaussian process regression (GPR)-based attack localization algorithm, from \cite{Zheng2022ResilientPrior}, to generate the support prior. The prior generator is to use GPR models to learn the nominal maps from hidden auxiliary states to the measurements. Attacked measurements are revealed if they cannot be explained by the trained GPR models with high likelihood. As discussed in Remark~\ref{rmk:compare_error_bounds} and shown in the numerical simulation section, better resiliency through the proposed weighted $\ell_1$ estimator is due to the precision guarantee of prior over $50\%$. However, given the probabilistic nature of data-driven approaches, relying solely on the prior derived from data-driven attack detectors may not rigorously ensure performance consistently. Therefore, we adopted a robustness mechanism, called pruning algorithm \cite{zheng2021Resilient}, to reinforce the prior. In this simulation, the mean of localization precision achieved $E(\textsf{PPV}) = 65.5\%$, as shown in Fig.~\ref{fig:prior_precision}. So the generated prior information is better than the random flip of a fair coin. As shown in Fig. \ref{Resilient_pruning_simulation}, we call the observer design, composed of the GPR-based localization algorithm, the pruning algorithm, and the weighted $\ell_1$ estimator, the resilient pruning observer (RPO).

 In this simulation, we compared RPO with other typical resilient observer designs, such as an event-triggered Luenberger observer (ETLO) \cite{shoukry2015event}, an unconstrained $\ell_1$ observer (UL1O) \cite{fawzi2014secure}, and a multi-model observer (MMO) \cite{anubi2019enhanced}. UL1O utilized the $\ell_1$ minimization program to recover the sparse attack signal \cite{fawzi2014secure}. ETLO is a resilient moving-horizon Luenberger-like observer solved by a projected gradient descent technique \cite{shoukry2015event}. MMO is a constrained $\ell_1$ estimator whose constraints include the system update law and the GPR attack detector \cite{anubi2019enhanced}. UL1O, ETLO, and MMO assume $2Tk$-observability which is rephrased in Corollary~\ref{Thm:recovery_Withoutprior} as a Row-RIP condition of order $(a+1)Tk$. This condition is satisfied when half of the sensors are safe \cite{fawzi2014secure,chong2015observability}. In this application example, $4$ of $19$ sensors are attacked. As further indicated in Theorem~\ref{Thm:recovery_Withprior}, the same-level Row-RIP condition, as the one in Corollary~\ref{Thm:recovery_Withoutprior}, is designed for our proposed weighted $\ell_1$ estimator. Thus, the condition of the proposed RPO is also satisfied.

Fig.~\ref{Comparison_all} and Table \ref{table_of_error_metric} show the estimation errors of generator rotor angles by $4$ resilient observers. LO becomes unstable in the presence of attacks. UTLO, ETLO, and MMO have non-negligible errors in terms of the magnitude of nominal states ($|\delta|<0.005)$. ETLO has smooth estimation errors since it is a dynamic observer design. MMO is an observer design that incorporates a prior measurement, so it works better than UL1O. However, the prior performance is not stable due to the uncertainty of the data-driven detector. By adding the pruning algorithm, the precision of prior information ($E(\textsf{PPV}) = 65.5\%$) is guaranteed to be better than a random flip of a fair coin. Therefore, RPO produces better estimation performance, which validates the statements in Remark~\ref{rmk:compare_error_bounds}. 
\begin{table} 
\centering
 \caption{Error Metric Values\label{table_of_error_metric}}
 \resizebox{0.7\textwidth}{!}{\begin{tabular}{|c | c | c | c | c | c | c|}
\hline
 &\multicolumn{5}{|c|}{RMS Metric} \\
 \hline
 & \textbf{LO} & \textbf{UL1O} & \textbf{MMO} & \textbf{ETLO} & \textbf{RPO}\\
 \hline
 $\delta_1$ & $2.953$ & $6.2e-5$ & $8.28e-6$ & $5.61e-4$ & $1.13e-15$ \\ 

$\delta_2$ & $2.881$ & $6.66e-5$ & $7.89e-6$ & $8.24e-4$ & $2.55e-16$ \\

 $\delta_3$ & $3.090$ & $5.45e-4$ & $4.65e-5$ & $6.06e-4$ & $8.98e-16$ \\

 $\delta_4$ & $3.195$ & $2.55e-4$ & $2.01e-5$ & $3.6e-3$ & $9.44e-16$ \\
 
 $\delta_5$ & $3.412$ & $6.41e-4$ & $1.95e-4$ & $9.66e-4$ & $6.62e-16$ \\  
\hline \hline
&\multicolumn{5}{|c|}{Max. Ans. Metric} \\
 \hline
 & \textbf{LO}  & \textbf{UL1O} & \textbf{MMO} & \textbf{ETLO} & \textbf{RPO}\\
 \hline
 $\delta_1$ & $9.729$ & $0.0017$ & $1.82e-4$ & $0.0012$ & $3.1e-14$ \\ 

$\delta_2$ & $9.482$ & $0.0017$ & $1.53e-4$ & $0.0017$ & $5.94e-15$  \\

 $\delta_3$ & $13.423$ & $0.0116$ & $8.37e-4$ & $0.0013$ & $2.45e-14$  \\

 $\delta_4$ & $12.834$ & $0.0058$ & $3.55e-4$ & $0.0079$ & $2.42e-14$ \\
 
 $\delta_5$ & $12.592$ & $0.0078$ & $0.0027$ & $0.0021$ & $1.28e-14$  \\  
 \hline
\end{tabular}}
\end{table}